\documentclass[preprint,10pt,3p]{elsarticle}

\usepackage{amssymb}
\usepackage{latexsym}
\usepackage{amsmath}
\usepackage{amsfonts}
\usepackage{amsthm}
\usepackage{mathrsfs}
\usepackage{bm}
\usepackage{graphicx}
\usepackage{subfigure}
\usepackage{float}
\usepackage{multirow}
\usepackage{color}

\newcommand{\Proj}{{\bm{P}}}
\newcommand{\cProj}{{\mathcal{P}}}
\newcommand{\cProjLH}{{\mathscr{P}}}
\newcommand{\dif}{\mathrm{d}}
\newcommand{\dt}{{\Delta t}}
\newcommand{\exOp}{{\bm{X}}^{[\text{E}]}}
\newcommand{\imOp}{{\bm{X}}^{[\text{I}]}}
\newcommand{\nStages}{n_{\tiny \textsf{s}}}
\newcommand{\Dim}{\scriptsize \textsc{D}}

\newcommand{\taubold}{{\boldsymbol \tau}}

\newcommand{\revise}{}

\journal{International Journal for Numerical Methods in Engineering}

\begin{document}

\begin{frontmatter}
  
\title{Fourth- and higher-order finite element methods
  for the incompressible Navier-Stokes equations
  with Dirichlet boundary conditions}

\author[zju]{Yang Li}
\ead{3150104952@zju.edu.cn}

\author[zju]{Heyu Wang}
\ead{wangheyu@zju.edu.cn}

\author[zju,future]{Qinghai Zhang\corref{corres}}
\ead{qinghai@zju.edu.cn}

\cortext[corres]{Corresponding author.}

\address[zju]{School of Mathematical Sciences,
  Zhejiang University,
  Hangzhou,
  Zhejiang,
  310058,
  China}

\address[future]{Institute of Fundamental and Transdisciplinary Research,
  Zhejiang University,
  Hangzhou,
  Zhejiang,
  310058,
  China}

\begin{abstract}
Inspired by the unconstrained pressure Poisson equation (PPE) formulation
[Liu, Liu, \& Pego, Comm. Pure Appl. Math. 60 (2007): 1443-1487],
we previously proposed
the generic projection and unconstrained PPE (GePUP) formulation
[Zhang, J. Sci. Comput. 67 (2016): 1134-1180]
for numerically solving the incompressible Navier-Stokes equations (INSE)
with no-slip boundary conditions.
In GePUP,
the main evolutionary variable does not have to be solenoidal
with its divergence controlled by a heat equation.
This work presents high-order finite-element solvers for the INSE
under the framework of method-of-lines.
Continuous Lagrange finite elements of equal order are utilized
for the velocity and pressure finite element spaces to
discretize the weak form of GePUP in space,
while high-order implicit-explicit Runge-Kutta methods are then employed to
treat the stiff diffusion term implicitly and
the other terms explicitly.
Due to the implicit treatment of the diffusion term,
the time step size is only restricted by convection.
The solver is efficient in that
advancing the solution at each time step only involves
solving a sequence of linear systems
either on the velocity or on the pressure
with geometric multigrid methods.
Furthermore,
the solver is enhanced with adaptive mesh refinement
so that the multiple length scales and time scales
in flows at moderate or high Reynolds numbers
can be efficiently resolved.
Numerical tests with various Reynolds numbers
are performed for
the single-vortex test,
\revise{the lid-driven cavity,
  and the flow past a cylinder/sphere,}
demonstrating the high-order accuracy of GePUP-FEM both in time and in space
and its capability of accurately and efficiently capturing the right physics.
Moreover,
our solver offers the flexibility in
choosing velocity and pressure finite element spaces
and is free of the standard inf-sup condition.
\end{abstract}


\end{frontmatter}


\section{Introduction}
\label{sec:introduction}

The Navier-Stokes equations (NSE) represent
a fundamental set of equations
governing the behavior of fluid flows and
play a pivotal role in understanding a wide range of phenomena
across various fields such as
engineering, physics, meteorology, and biology.
These equations,
derived from the conservation laws of mass, momentum, and energy,
describe the motion of viscous fluids
under the influence of external forces.
Understanding and solving the NSE are crucial for
predicting the behavior of fluids in diverse scenarios,
ranging from the airflow over an aircraft wing to
the blood flow in human arteries.
Moreover,
these equations serve as the foundation for simulating and optimizing
complex fluid dynamics problems
encountered in practical engineering applications,
including aerodynamics, hydrodynamics, heat transfer,
and chemical engineering processes.
Therefore,
the development and application of efficient and accurate numerical methods
for solving the NSE
not only advance our understanding of fluid mechanics,
but also enable the design and optimization of innovative technologies
that shape modern industries and enhance our everyday lives.

The incompressible Navier-Stokes equations (INSE)
with Dirichlet boundary conditions read
\begin{subequations}
  \label{eq:INSE}
  \begin{alignat}{2}
    \frac{\partial \bm{u}}{\partial t} + \bm{u}\cdot\nabla\bm{u}
    &= \bm{f} - \nabla p + \nu\Delta\bm{u} &\quad& \text{in } \Omega, \\
    \nabla\cdot\bm{u} &= 0 && \text{in } \Omega, \\
    \bm{u} &= \bm{g} && \text{on } \partial\Omega, \\
    \bm{u}(\bm{x}, 0) &= \bm{u}_0(\bm{x}) && \text{in } \Omega,
  \end{alignat}
\end{subequations}
where $t$ is time,
$\bm{x}\in \Omega\subset \mathbb{R}^{\Dim} \, (\Dim=2,3)$
is the spatial location,
$\bm{f}$ is the external forcing term,
$\bm{g}$ is the given velocity boundary function
satisfying $\int_{\partial \Omega}\bm{g}\cdot\bm{n} = 0$,
$\bm{u}_0$ is the given initial velocity field
assumed to be divergence-free,
$\bm{u}$ is the velocity,
$p$ is the pressure,
and $\nu$ is the kinematic viscosity.

As a fundamental and notoriously difficult problem,
the well-posedness of the INSE
(in $\mathbb{R}^3$ and $\mathbb{R}^3/\mathbb{Z}^3$
\cite{fefferman06:_exist_navier_stokes})
is selected by Smale \cite{smale98:_mathem}
on his list of 18 mathematical problems for the 21st century
and by the Clay Institute
on the list of 7 millennium problems \cite{devlin03}.
Since no one has solved this problem yet,
numerical computation appears to be the only viable way for now
to obtain approximate solutions of the INSE.

\subsection{Challenges for the design of numerical solvers}
\label{sec:chall-design-numer}

In developing numerical methods for solving the INSE,
we are confronted with four major challenges.
\begin{enumerate}[(A)]
\item
  How to fulfill the divergence-free constraint (\ref{eq:INSE}b)?
\item
  How to achieve high-order accuracy both in time and in space?
\item
  How to accurately and efficiently capture
  structures of multiple length scales and time scales
  typically developed in flows with moderate or high Reynolds numbers?
\item
  How to decouple temporal integration from spatial discretization so that
  (i) the entire solver is constituted by orthogonal modules for these aspects,
  and
  (ii) solution methods for each aspect can be employed
  in a black-box manner
  and thus easily changed to make the entire solver versatile?
\end{enumerate}

Challenge (A) concerns
the special role of the pressure as a Lagrangian multiplier
in enforcing the divergence-free constraint (\ref{eq:INSE}b).
The lack of an evolution equation and suitable boundary conditions
of the pressure makes it difficult to obtain
high-order solutions for the pressure.

Challenges (B) and (C) concern accuracy and efficiency.
A numerical method should resolve all scales
that are relevant to the important physics.
Compared with fourth- and higher-order methods,
first- and second-order methods have simpler algorithms and
cheaper computations,
but towards a given accuracy
the computational resources may be rapidly exhausted.
It is shown both theoretically and analytically
in \cite[Section 7]{zhang16:_GePUP} that
fourth-order methods may have a significant efficiency advantage
over second-order methods.
In fact,
the speedup of high-order methods over low-order methods,
as measured by the ratio of their CPU times, 
grows as a power function
when the grid size or the targeted accuracy is reduced.

Challenge (B) also concerns faithful simulations of flows
where velocity derivatives such as vorticity
crucially affect the physics.
For first-order methods,
the computed velocity converges,
but the vorticity does not,
nor does the velocity gradient tensor.
Consequently,
the $O(1)$ error in the velocity gradient may
lead to structures different from
that of the original flow.
Thus it is not clear
whether solutions of a first-order method converge to
the \emph{right} physics.
Similar suspicions apply to second-order methods for flows
where second derivatives of the velocity are important.
In particular,
high-order accurate solvers are essential
for simulating fluid-structure interaction problems, 
where accurate tractions are needed at fluid-solid interfaces.

Challenge (D) concerns versatility and user-friendliness
of the numerical solver.
To cater for the problem at hand,
it is often desirable to change the time integrator from one to another.
For example,
flows with high viscosity are usually stiff
while those with small viscosity are not;
accordingly,
an implicit time integrator should be used in the former case
while an explicit one is usually suitable for the latter.
If the internal details of a time integrator are rigidly coupled into the INSE solver,
changing the time integrator would be
difficult and highly inconvenient.
Hence a time integrator should be treated as a black box:
for the ordinary differential equation (ODE)
$\frac{\dif \bm{U}}{\dif t} = \bm{f}(\bm{U}, t)$,
we should only need to feed into the time integrator
the initial condition $\bm{U}^n$ and
samples of $\bm{f}$ at a number of time instances to
get the solution $\bm{U}^{n+1}$ from the black box.
The versatility further leads to user-friendliness.
Analogous to an orthogonal basis of a vector space,
the mutually independent policies span a space of solvers.
In other words,
each solver can be conveniently assembled by
selecting a module for each constituting policy.

\subsection{Some previous methods}
\label{sec:some-prev-meth}

\subsubsection{Mixed finite element methods}
\label{sec:mixed-finite-element}

Mixed finite element methods
\cite{boffi13:_mixed_finit_elemen_method_applic}
tackle challenge (A) by treating
both the velocity and the pressure simultaneously
in an implicit fashion.
For example,
using the backward Euler method for time discretization and
neglecting the nonlinear convection term,
we are tasked with solving the following saddle-point problem
at each time step:
find $\bm{u}_h^{n+1}\in \bm{V}_h$ and $p_h^{n+1}\in Q_h$
such that
\begin{subequations}
  \label{eq:mixedFEM}
  \begin{alignat}{2}
    \forall \bm{v}_h\in \bm{V}_h,
    &\quad&
    \frac{1}{\dt}(\bm{u}_h^{n+1}, \bm{v}_h)_{\Omega}
    + \nu(\nabla\bm{u}_h^{n+1}, \nabla\bm{v}_h)_{\Omega}
    - (p_h^{n+1}, \nabla\cdot\bm{v}_h)_{\Omega} &=
    \frac{1}{\dt}(\bm{u}_h^n, \bm{v}_h)_{\Omega}
    + (\bm{f}^{n+1}, \bm{v}_h)_{\Omega}, \\
    \forall q_h\in Q_h,
    &\quad&
    (q_h, \nabla\cdot\bm{u}_h^{n+1})_{\Omega} &= 0,
  \end{alignat}
\end{subequations}
where $\bm{V}_h$ and $Q_h$ are the finite element spaces
for the velocity and the pressure,
respectively,
$\dt$ is the time step size,
and $(\cdot, \cdot)_{\Omega}$ denotes
the standard scalar-, vector-, or tensor-valued $L^2$ inner product
defined on $\Omega$.

Over the last fifty years,
the mixed finite element methods have achieved tremendous successes
and established elegant mathematical theories,
the most significant of which is the inf-sup condition,
also known as the Ladyzhenskaya-Babu\v{s}ka-Brezzi (LBB) condition,
\begin{equation}
  \label{eq:infSupCondition}
  \exists\beta > 0 \text{ s.t. }
  \forall h,
  \inf_{q_h\in Q_h}\sup_{\bm{v}_h\in \bm{V}_h}
  \frac{(q_h, \nabla\cdot\bm{v}_h)_{\Omega}}
  {\|q_h\|_{L^2(\Omega)}\|\bm{v}_h\|_{\bm{H}^1(\Omega)}} \ge \beta,
\end{equation}
which is both sufficient and necessary for
the well-posedness of the saddle-point problem (\ref{eq:mixedFEM}).
However,
this condition limits the flexibility of the numerical methods
as it is well known that
equal-order Lagrange finite-element pairs
for velocity and pressure approximations
violate the inf-sup condition (\ref{eq:infSupCondition}).
Although it is possible to
circumvent the inf-sup condition (\ref{eq:infSupCondition})
by employing some stabilization techniques
(such as continuous interior penalty
\cite{burman07:_contin_navier_stokes}
or subgrid-scale \cite{codina07:_time})
at the cost of additional complexity,
developing efficient solution techniques
for the saddle-point system (\ref{eq:mixedFEM})
still poses significant challenges
\cite{benzi05:_numer}.

\subsubsection{Projection methods}
\label{sec:projection-methods}

In stark contrast to the mixed finite element methods,
the projection methods,
pioneered by Chorin \cite{chorin68:_numer_solut_navier_stokes_equat} and
Temam \cite{temam69:_sur_navier_stokes},
stand out for their high efficiency
among various numerical methods for the INSE;
they are based on the observation that
the momentum equation (\ref{eq:INSE}a) can be rewritten as
\begin{equation}
  \label{eq:INSEa-EulerianAccel}
  \bm{a}^* = \bm{a} + \nabla p = \cProjLH\bm{a}^* + \nabla p,
\end{equation}
where $\bm{a}$ and $\bm{a}^*$ are the \emph{Eulerian acceleration}s
defined by
\begin{equation}
  \label{eq:EulerianAccelerations}
  \bm{a} := \frac{\partial \bm{u}}{\partial t},\qquad
  \bm{a}^* := -\bm{u}\cdot\nabla\bm{u}
  +\bm{f} +  \nu\Delta \bm{u},
\end{equation}
and $\cProjLH$ is the Leray-Helmholtz projection operator defined by
\begin{equation*}
  \bm{v}^* = \cProjLH\bm{v}^* + \nabla\phi := \bm{v} + \nabla\phi,
\end{equation*}
where
$\nabla\cdot\bm{v} = 0$ in $\Omega$ and
$\bm{v}\cdot\bm{n}$ is given on $\partial\Omega$
and satisfies $\int_{\partial\Omega}\bm{v}\cdot\bm{n} = 0$.

In the original projection method independently proposed by
Chorin \cite{chorin68:_numer_solut_navier_stokes_equat}
and Temam \cite{temam69:_sur_navier_stokes},
the initial condition $\bm{u}^n\approx \bm{u}(t^n)$
is first advanced to an auxiliary velocity $\bm{u}^*$
without worrying about the pressure gradient term
and then $\bm{u}^*$ is projected to the divergence-free space to
obtain $\bm{u}^{n+1}$,
\begin{subequations}
  \label{eq:1stOrderProjection}
  \begin{align}
    \frac{\bm{u}^*-\bm{u}^n}{\dt}
    &= -\bm{C}(\bm{u}^*, \bm{u}^n)
      + \bm{f}^n + \nu\bm{L}\bm{u}^*, \\
    \bm{u}^{n+1} &= \Proj\bm{u}^*,
  \end{align}
\end{subequations}
where $t^n$ is the starting time of the $n$th step,
\mbox{$\bm{f}^n= \bm{f}(t^n)$},
$\bm{C}(\bm{u}^*, \bm{u}^n)\approx
[\bm{u}\cdot\nabla\bm{u}](t^n)$,
and
$\bm{L}$ and $\Proj$ 
are discrete approximations of the Laplacian $\Delta$
and the Leray-Helmholtz projection $\cProjLH$, respectively.
As a result,
the INSE is advanced by
solving a sequence of elliptic boundary value problems (BVPs)
either on the velocity or on the pressure.

Chorin and Temam's original projection method (\ref{eq:1stOrderProjection}) is
only first-order accurate in time and
its improvement to the second order has been
the aim of many subsequent works,
see the two excellent reviews
\cite{brown01:_accurate_projection_methods_for_ins,
  guermond06:_overview_of_projection_methods_for_incompressible_flows}.
Moreover,
there also exist many second-order adaptive projection methods
in the literature,
e.g.,
\cite{howell97,
  martin08:_navier_stokes,
  trebotich15:_navier_stokes,
  blomquist24:_stabl}.

A common basis of many second-order methods
is the temporal discretization of the INSE (\ref{eq:INSE}) with the trapezoidal rule,
\begin{subequations}
  \label{eq:2ndOrderProjectionIdea}
  \begin{align}
    \frac{\bm{u}^{n+1}-\bm{u}^n}{\dt} + \nabla p^{n+\frac{1}{2}}
    &= -[\bm{u}\cdot\nabla\bm{u}]^{n+\frac{1}{2}}
      + \bm{f}^{n+\frac{1}{2}}
      + \frac{\nu}{2}\Delta(\bm{u}^{n+1}+\bm{u}^n), \\
    \nabla\cdot \bm{u}^{n+1} &= 0,
  \end{align}
\end{subequations}
where  $p^{n+\frac{1}{2}}\approx p(t^{n+\frac{1}{2}})$,
$\bm{f}^{n+\frac{1}{2}} = \bm{f}(t^{n+\frac{1}{2}})$,
and $[\bm{u}\cdot\nabla\bm{u}]^{n+\frac{1}{2}}\approx
[\bm{u}\cdot\nabla\bm{u}](t^{n+\frac{1}{2}})$
are numerical approximations
at the mid-time $t^{n+\frac{1}{2}}=\frac{1}{2}(t^n+t^{n+1})$
of the $n$th step.
Based on (\ref{eq:2ndOrderProjectionIdea}),
many second-order projection methods
with fractional stepping can be united 
\cite{brown01:_accurate_projection_methods_for_ins} as follows,
\begin{itemize}
\item
  choose an auxiliary variable $q$ and
  a boundary condition $\bm{B}(\bm{u}^*) = \bm{0}$
  to solve for $\bm{u}^*$,
  \begin{equation}
    \label{eq:intermediateFieldSolve}
    \frac{\bm{u}^*-\bm{u}^n}{\dt} + \nabla q
    = -[\bm{u}\cdot\nabla\bm{u}]^{n+\frac{1}{2}}
    + \bm{f}^{n+\frac{1}{2}}
    + \frac{\nu}{2}\Delta(\bm{u}^*+\bm{u}^n);
  \end{equation}
\item
  enforce (\ref{eq:2ndOrderProjectionIdea}b)
  with a discrete projection so that
  \begin{equation}
    \label{eq:2ndOrderProjectionStep}
    \bm{u}^* = \bm{u}^{n+1} + \dt \nabla\phi^{n+1};
  \end{equation}
\item
  update the pressure by
  \begin{equation}
    \label{eq:pressureUpdate}
    p^{n+\frac{1}{2}} = q + U(\phi^{n+1}),
  \end{equation}
\end{itemize}
where the specific forms of $q$, $\bm{B}(\bm{u}^*)$,
and $U(\phi^{n+1})$
vary from one method to another.
Two examples are the method of Bell, Colella, and Glaz
\cite{bell89:_navier_stokes}
specified by
\begin{equation}
  \label{eq:BellColellaGlaz}
  q = p^{n-\frac{1}{2}}, \quad
  \bm{u}^*|_{\partial\Omega} = \bm{g}^{n+1},
  \quad
  \nabla p^{n+\frac{1}{2}} = \nabla p^{n-\frac{1}{2}} + \nabla\phi^{n+1},
\end{equation}
and that of Kim and Moin
\cite{kim85:_applic_navier_stokes}
by
\begin{subequations}
  \label{eq:KimMoin}
  \begin{align}
    q = 0, \quad
    \bm{n}\cdot\bm{u}^*|_{\partial\Omega}
    &= \bm{n}\cdot\bm{g}^{n+1}, \quad
      \taubold\cdot\bm{u}^*|_{\partial\Omega}
      = \taubold\cdot(\bm{g}^{n+1} + \dt\nabla\phi^n)|_{\partial\Omega}, 
    \\
    \nabla p^{n+\frac{1}{2}}
    &= \nabla\phi^{n+1} - \frac{\nu \dt}{2}\nabla\Delta\phi^{n+1}.
  \end{align}
\end{subequations}

As observed by Brown, Cortez, and Minion
\cite{brown01:_accurate_projection_methods_for_ins},
the choices of $q$, $\bm{B}(\bm{u}^*)$,
and $U(\phi^{n+1})$ are not independent.
Indeed,
substitute (\ref{eq:2ndOrderProjectionStep})
into (\ref{eq:intermediateFieldSolve}),
subtract (\ref{eq:2ndOrderProjectionIdea}a),
and we have
\begin{equation}
  \label{eq:exactPressureUpdate}
  \nabla p^{n+\frac{1}{2}} = \nabla q +\nabla\phi^{n+1}
  -\frac{\nu \dt}{2}\Delta\nabla\phi^{n+1}.
\end{equation}
Consequently,
the pressure computed by (\ref{eq:BellColellaGlaz})
is only first-order accurate
\cite{brown01:_accurate_projection_methods_for_ins}
because its pressure update formula
is at an $O(\Delta t)$ discrepancy with (\ref{eq:exactPressureUpdate}).
On the other hand,
the pressure update formula (\ref{eq:KimMoin}b)
complies with (\ref{eq:exactPressureUpdate})
due to the commutativity of the Laplacian and the gradient,
but (\ref{eq:KimMoin}b) and
(\ref{eq:2ndOrderProjectionStep})
imply $ p^{n+\frac{1}{2}} = \phi^{n+1} - \frac{\nu}{2}
\nabla\cdot\bm{u}^*$,
which may deteriorate the pressure accuracy to the first order
\cite{brown01:_accurate_projection_methods_for_ins}.
Nonetheless,
both methods (\ref{eq:BellColellaGlaz}) and (\ref{eq:KimMoin})
are second-order accurate for the velocity.

Projection methods have been successful.
However,
enhancing them to fourth- and higher-order accuracy
presents a challenge.
As discussed above, 
the choices of $q$, $\bm{B}(\bm{u}^*)$, and $U(\phi^{n+1})$
are coupled according to internal details of the time integrator.  
Consequently,
switching from one time integrator to another
calls for a new derivation.
Furthermore,
although appearing divorced,
the velocity and the pressure are still implicitly coupled
by the boundary condition of
the auxiliary variable $\bm{u}^*$,
with the coupling determined not by physics
but still by internal details of the time integrator.
Therefore, 
it is highly challenging, if not impossible,
to answer challenges (B) and (D)
via the approach of traditional projection methods.

\subsection{The formulation of pressure Poisson equation (PPE)}
\label{sec:form-press-poiss}

The PPE describes an \emph{instantaneous} relation
between the pressure and the velocity in the INSE
and, in domains with Dirichlet boundary conditions, has the form
\begin{subequations}
  \label{eq:PPE}
  \begin{alignat}{2}
    \Delta p &= \nabla \cdot
    (\bm{f}-\bm{u}\cdot\nabla\bm{u}
    + \nu \Delta \bm{u})
    &\quad& \text{in } \Omega,
    \\
    \bm{n}\cdot \nabla p
    &= \bm{n}\cdot
    \left(\bm{f} - \bm{u}\cdot\nabla\bm{u} + \nu\Delta\bm{u}
      - \frac{\partial \bm{g}}{\partial t}\right)
    && \text{on } \partial \Omega,
  \end{alignat}
\end{subequations}
where (\ref{eq:PPE}a) follows from the divergence of (\ref{eq:INSE}a) and
the divergence-free condition (\ref{eq:INSE}b),
while (\ref{eq:PPE}b) from
the normal component of (\ref{eq:INSE}a) on $\partial\Omega$ and
the Dirichlet boundary condition (\ref{eq:INSE}c).
For the PPE in domains with other boundary conditions,
(\ref{eq:PPE}b) should be replaced with
the normal component of (\ref{eq:INSEa-EulerianAccel}).

(\ref{eq:INSE}a), (\ref{eq:INSE}c),
(\ref{eq:PPE}), and the additional boundary condition
$\nabla\cdot\bm{u}=0$ on $\partial \Omega$
are collectively called
\emph{the PPE formulation of the INSE
  with Dirichlet boundary conditions},
which is equivalent to the original INSE
\cite{gresho87:_navier_stokes}.
In terms of computation, however, 
the PPE formulation has a decisive advantage over the original INSE.
If (\ref{eq:INSE}a) is discretized in time
with (\ref{eq:INSE}b) as a constraint,
the resulting index-2 differential algebraic system
may suffer from large order reductions
\cite{sanderse12:_accur_runge_kutta_navier_stokes}.
In contrast,
replacing the divergence-free constraint
with the PPE avoids this difficulty.
As such,
the PPE formulation allows
the time integrator to be treated as a black box
and thus to be easily changed;
indeed, the pressure is an implicit function of $\bm{u}$
and its interaction with $\bm{u}$
is completely decoupled from
internal details of the time integrator.
Also,
there is no need to introduce
nonphysical auxiliary variables.
\revise{
These advantages of the PPE formulation have led to
the development of numerous successful numerical methods.
Two prominent examples are
the unconstrained PPE formulation
\cite{johnston04:_accur_navier_stokes,
  liu07:_stabil_conver_effic_navier_stokes,
  liu09:_error_navier_stokes,
  liu09:_open_navier_stokes,
  liu10:_stable_accurate_pressure_unsteady_incompressible_viscous_flow,
  jia11:_stabl_navier_stokes}
and the PPE formulation with electric boundary conditions
\cite{shirokoff11:_navier_stokes,
  rosales21:_high_poiss_navier_stokes},
which we briefly review in Sections \ref{sec:form-unconstr-ppe} and
\ref{sec:ppe-formulation-with},
respectively.
Earlier contributions include:
Karniadakis et al. \cite{karniadakis91:_high_navier_stokes},
who solved the PPE (\ref{eq:PPE}) by
replacing the diffusion term $\nu\Delta\bm{u}$ in (\ref{eq:PPE}b)
with $-\nu\nabla\times\nabla\times\bm{u}$ and
employing a temporal extrapolation of the velocity field;
Petersson \cite{petersson01:_stabil_stokes_navier_stokes},
who investigated the stability properties of different pressure boundary conditions;
and Leriche et al. \cite{leriche06:_numer_stokes},
who proposed a non-splitting scheme by
introducing an intermediate divergence-free acceleration field
$\bar{\bm{a}} = \partial \bm{u}/\partial t - \nu\Delta\bm{u}$.
More recent contributions include:
Li \cite{li20:_navier_stokes},
who investigated accurate pressure boundary conditions at the discrete level;
Meng et al. \cite{meng20:_fourt_imex_navier_stokes},
who addressed complex moving domain problems
using a WENO-based scheme on composite overlapping grids;
and Pacheco et al. \cite{pacheco21:_newton},
who focused on simulating non-Newtonian incompressible flow problems.
}

Unfortunately,
as observed by Liu, Liu, and Pego
\cite{liu07:_stabil_conver_effic_navier_stokes},
(\ref{eq:INSE}a) and (\ref{eq:PPE}a)
yield
\begin{equation}
  \label{eq:degenerateVelDivPPE}
  \frac{\partial \nabla\cdot\bm{u}}{\partial t} = 0;
\end{equation}
this degenerate equation implies that
in the PPE formulation
we have no control over $\nabla\cdot\bm{u}$
and its evolution is up to the particularities
of the numerical schemes.

\subsection{The formulation of unconstrained PPE (UPPE)}
\label{sec:form-unconstr-ppe}

By applying the Leray-Helmholtz projection $\cProjLH$
to (\ref{eq:INSEa-EulerianAccel}),
Liu, Liu, and Pego
\cite{liu07:_stabil_conver_effic_navier_stokes}
obtained
\begin{equation}
  \label{eq:UPPE1}
  \frac{\partial \bm{u}}{\partial t}
  - \cProjLH \bm{a}^*
  = \nu\nabla(\nabla\cdot \bm{u}),
\end{equation}
where the zero right-hand side (RHS) is 
added for stability reasons. 
The divergence of (\ref{eq:UPPE1}) gives
\begin{equation}
  \label{eq:heatEqInUPPE}
  \frac{\partial(\nabla\cdot\bm{u})}{\partial t}
  = \nu\Delta(\nabla\cdot\bm{u}),
\end{equation}
which, together with the boundary condition
$\nabla\cdot \bm{u}=0$ on $\partial\Omega$,
dictates an exponential decay of 
a nonzero $\nabla\cdot \bm{u}$ in $\Omega$.
A juxtaposition of (\ref{eq:heatEqInUPPE})
and (\ref{eq:degenerateVelDivPPE})
exposes a prominent advantage of 
(\ref{eq:UPPE1})
that any divergence residue caused by discretization errors
is now well over control.

With the identity $\nabla(\nabla\cdot \bm{u})
=\Delta({\mathcal I}-\cProjLH)\bm{u}$
and the Laplace-Leray commutator
$[\Delta, \cProjLH] := \Delta \cProjLH - \cProjLH \Delta$,
Liu, Liu, and Pego
\cite{liu07:_stabil_conver_effic_navier_stokes}
rewrote (\ref{eq:UPPE1}) as 
\begin{equation}
  \label{eq:UPPE2}
  \frac{\partial \bm{u}}{\partial t}
  + \cProjLH (\bm{u}\cdot\nabla \bm{u} -\bm{f})
  + \nu [\Delta, \cProjLH]\bm{u}
  = \nu\Delta \bm{u}, 
\end{equation}
which provides a novel perspective on the INSE
as a controlled perturbation of the vector diffusion equation
$\frac{\partial \bm{u}}{\partial t} = \nu\Delta \bm{u}$.
For a bounded and connected domain $\Omega$
with ${\mathcal C}^3$ boundary,
they established a sharp estimate on
the $L^2$-norm of the Laplace-Leray commutator $[\Delta, \cProjLH]\bm{u}$,
\begin{equation}
  \label{eq:estimateOnLaplaceLeray}
  \forall\varepsilon>0,
  \exists C \ge 0
  \text{ s.t. }
  \forall \bm{u}\in \bm{H}^2(\Omega)\cap\bm{H}^1_0(\Omega), \, \,
  \int_{\Omega}\left| [\Delta, \cProjLH]\bm{u} \right|^2 \le
  \left( \frac{1}{2}+\varepsilon \right)\int_{\Omega}\left| \Delta\bm{u} \right|^2
  + C\int_{\Omega}|\nabla\bm{u}|^2
\end{equation}
and proved the unconditional stability and convergence
of a first-order scheme
without the standard inf-sup condition (\ref{eq:infSupCondition}),
\begin{subequations}
  \label{eq:UPPE-1stOrder}
  \begin{align}
    \left(\nabla p^n, \nabla \eta\right)_{\Omega}
    &= \left(\bm{f}^{n}- \bm{u}^n\cdot\nabla \bm{u}^{n}
      +\nu\Delta \bm{u}^{n} - \nu\nabla\nabla\cdot\bm{u}^{n},
      \nabla \eta\right)_{\Omega},
    \\
    \frac{\bm{u}^{n+1}-\bm{u}^n}{\dt} - \nu\Delta\bm{u}^{n+1}
    &= \bm{f}^{n}- \bm{u}^n\cdot\nabla \bm{u}^{n}
      -\nabla p^n \qquad \text{ in } \Omega,\\
    \bm{u}^{n+1} &= \bm{0} \qquad\qquad\qquad\qquad\qquad\quad
                   \text{ on } \partial\Omega,
  \end{align}
\end{subequations}
where (\ref{eq:UPPE-1stOrder}a) is the PPE in weak form with
the test function
$\eta \in H^1_{\int=0}(\Omega) :=
\left\{ \eta\in H^1(\Omega): \int_{\Omega}\eta = 0 \right\}$.
For simplicity,
(\ref{eq:UPPE-1stOrder}) is presented for the case
where the velocity satisfies homogeneous boundary conditions; 
nonhomogeneous conditions are discussed
in \cite{liu07:_stabil_conver_effic_navier_stokes}.

In (\ref{eq:UPPE-1stOrder}),
the computations of pressure and velocity have been decoupled by
treating the pressure gradient term explicitly
and the viscous term implicitly.
Its stability follows from 
(\ref{eq:estimateOnLaplaceLeray})
and the fact that
the pressure gradient in the INSE consists of two parts:
\begin{equation*}
  \nabla p = \nabla p_c + \nu\nabla p_s, \quad
  \nabla p_c := (\mathcal{I}-\cProjLH)(\bm{f} - \bm{u}\cdot\nabla\bm{u}), \quad
  \nabla p_s := [\Delta, \cProjLH]\bm{u},
\end{equation*}
where $\nabla p_c$ balances
the divergence of the forcing term and the nonlinear convection term
while $\nabla p_s$ accounts for
the Laplace-Leray commutator.
The stability and error estimates of the resulting fully discrete scheme
with $C^1$ finite elements for velocity
are further proved in \cite{liu09:_error_navier_stokes},
where the pair of velocity and pressure finite element spaces
need not satisfy the inf-sup condition (\ref{eq:infSupCondition}).
This has also been observed for $C^0$ finite element schemes
in numerical experiments
\cite{johnston04:_accur_navier_stokes,
  liu10:_stable_accurate_pressure_unsteady_incompressible_viscous_flow,
  li20:_navier_stokes}.

From (\ref{eq:UPPE2}),
Liu, Liu, and Pego
\cite{liu10:_stable_accurate_pressure_unsteady_incompressible_viscous_flow}
also deduced a strong form of the UPPE formulation as
\begin{subequations}
  \label{eq:UPPEstrong}
  \begin{alignat}{2}
    \frac{\partial \bm{u}}{\partial t}
    + \bm{u}\cdot\nabla\bm{u}
    &= \bm{f} -\nabla p +  \nu\Delta \bm{u}
    &\quad& \text{in } \Omega,
    \\
    \bm{u} &= \bm{g} &&
    \text{on } \partial\Omega,
    \\
    \Delta p
    &= \nabla\cdot(\bm{f}-\bm{u}\cdot\nabla\bm{u})
    &&  \text{in } \Omega, \\
    \bm{n}\cdot\nabla p
    &= \bm{n}\cdot \left( \bm{f} - \bm{u}\cdot\nabla\bm{u}
      - \nu\nabla\times\nabla\times\bm{u}
      - \frac{\partial \bm{g}}{\partial t} \right)
    && \text{on } \partial\Omega.
  \end{alignat}
\end{subequations}
Based on this formulation,
they further devised improved approximations
for the computation of pressure
in existing pressure-approximation and pressure-update projection methods,
and developed a slip-corrected projection method
\cite{liu10:_stable_accurate_pressure_unsteady_incompressible_viscous_flow}
that is third-order accurate
in time for both velocity and pressure.

The PPE (\ref{eq:PPE})
and the UPPE (\ref{eq:UPPEstrong}c,d)
have slightly different forms
and nonetheless a crucial distinction:
with the divergence of (\ref{eq:UPPEstrong}a),
(\ref{eq:UPPEstrong}c) leads to (\ref{eq:heatEqInUPPE})
whereas
(\ref{eq:PPE}a) leads to (\ref{eq:degenerateVelDivPPE}).

However,
as discussed in
\cite{liu09:_error_navier_stokes,
  liu09:_open_navier_stokes,
  liu10:_stable_accurate_pressure_unsteady_incompressible_viscous_flow},
a straightforward finite-element discretization of
the UPPE formulation (\ref{eq:UPPEstrong})
can lead to qualitatively incorrect results
for problems where the solution is not sufficiently smooth,
such as those with re-entrant corners.
One possible remedy is to
employ the Leray-Helmholtz projection operator $\cProjLH$ to
suppress the velocity divergence.
However,
$\cProjLH$ is absent in (\ref{eq:UPPEstrong})
and thus any projection on the velocity
in a method-of-lines (MOL) algorithm
would be a mismatch of the numerical algorithm
to the governing equations.
Moreover,
it is not clear
which $\bm{u}$'s in (\ref{eq:UPPEstrong}) should be projected in MOL.

Another issue with the UPPE formulation (\ref{eq:UPPEstrong}) is that
it is difficult for a discrete projection $\Proj$
with fourth- or higher-order accuracy
to satisfy properties of the Leray-Helmholtz projection $\cProjLH$,
such as
\begin{equation*}
  \cProjLH^2 = \cProjLH, \quad
  \nabla\cdot\cProjLH\bm{v}^* = 0, \quad
  \cProjLH\nabla\phi = \bm{0}.
\end{equation*}
In particular,
the discretely projected velocity might not be divergence-free.
Then how does the approximation error of
$\Proj$ to $\cProjLH$
affect the stability of the system of ODEs
under the MOL framework?
It is neither clear nor trivial
how to answer this question with (\ref{eq:UPPEstrong}).

\subsection{The PPE formulation
  with electric boundary conditions}
\label{sec:ppe-formulation-with}

To endow the heat equation (\ref{eq:heatEqInUPPE})
with homogeneous Dirichlet boundary conditions,
Shirokoff and Rosales \cite{shirokoff11:_navier_stokes}
proposed a novel PPE reformulation of the INSE,
incorporating electric boundary conditions (EBC):
\begin{subequations}
  \label{eq:PPEElectricBC}
  \begin{align}
    \frac{\partial \bm{u}}{\partial t} + \bm{u}\cdot\nabla\bm{u}
    &= \bm{f} - \nabla p + \nu\Delta\bm{u} \qquad\,\,\, \text{ in } \Omega, \\
    \bm{n}\times\bm{u} &= \bm{n}\times\bm{g},
    \quad \nabla\cdot\bm{u} = 0 \quad \text{ on } \partial\Omega, \\
    \Delta p &=
    \nabla\cdot (\bm{f} - \bm{u}\cdot\nabla\bm{u}) \qquad \text{ in } \Omega, \\
    \bm{n}\cdot\nabla p &=
    \bm{n}\cdot \left( \bm{f} - \bm{u}\cdot\nabla\bm{u}
      - \nu\nabla\times\nabla\times\bm{u}
      - \frac{\partial \bm{g}}{\partial t} \right)
    + \lambda\bm{n}\cdot(\bm{u} - \bm{g})  \quad \text{ on } \partial\Omega,
  \end{align}
\end{subequations}
where (\ref{eq:PPEElectricBC}b) is the EBC
commonly encountered in electrostatics,
and the term $\lambda\bm{n}\cdot(\bm{u} - \bm{g})$,
with $\lambda$ being a positive relaxation parameter,
ensures that
the normal velocity on the boundary decays exponentially.

To capture the structure imposed by 
the EBC (\ref{eq:PPEElectricBC}b),
Rosales et al. \cite{rosales21:_high_poiss_navier_stokes}
introduced the vorticity $\bm{\omega} = \nabla\times\bm{u}$
as a new variable and
employed mixed finite elements
to solve for $\bm{\omega}$ and $\bm{u}$ in space.
However,
this approach leads to a saddle-point problem,
demanding that the finite element pair
$(\bm{\omega},\bm{u})$
satisfy the inf-sup condition.
Moreover,
this mixed formulation also
increases the degrees of freedom (DoFs)
in the discrete problem.

\subsection{The contribution of this work}
\label{sec:contr-this-work}

Inspired by the UPPE formulation (\ref{eq:UPPEstrong}),
we previously proposed
the generic projection and unconstrained PPE (GePUP) formulation
\cite{zhang16:_GePUP}
for solving the INSE with no-slip boundary conditions.
This work is a continuation of \cite{zhang16:_GePUP} and
 our first step towards
the development of high-order accurate finite-element solvers
for simulating viscous incompressible flows.
Specifically,
we propose GePUP-FEM,
INSE finite-element solvers
that address all the challenges listed in Section \ref{sec:chall-design-numer}.
Novel features of GePUP-FEM include
\begin{enumerate}[(i)]
\item
  \emph{versatility}:
  time integration and spatial discretization are completely decoupled,
  so that high-order Runge-Kutta methods can be employed in a black-box manner;
\item
  \emph{accuracy}:
  fourth- and higher-order accuracy both in time and in space is achieved;
\item
  \emph{efficiency}:
  advancing the solution at each time step only involves
  solving a sequence of linear systems either on the velocity or on the pressure
  with geometric multigrid methods;
  moreover,
  GePUP-FEM is enhanced with adaptive mesh refinement to address Challenge (C);
\item
  \emph{flexibility}:
  the choice of finite element spaces for the velocity and the pressure
  is free of the inf-sup condition (\ref{eq:infSupCondition}).
\end{enumerate}
To the best of our knowledge,
GePUP-FEM is the first numerical method to possess all these features.

The rest of this paper is organized as follows.
In Section \ref{sec:gepup-formulation},
we briefly review the GePUP formulation of the INSE.
In Section \ref{sec:numerical-algorithms},
we derive the weak form of the GePUP formulation
and design high-order GePUP-FEM schemes
via the method of lines.
Benchmark problems are numerically solved
in Section \ref{sec:numerical-tests}
to confirm the high-order accuracy of GePUP-FEM and
its capability in accurately and efficiently capturing the right physics.
Finally,
Section \ref{sec:conclusion} concludes this paper
with some future research prospects.

\section{The GePUP formulation}
\label{sec:gepup-formulation}

In numerically simulating incompressible flows,
the computed velocity field rarely satisfies
the divergence-free condition \emph{exactly}.
Correspondingly,
we might as well
relax the evolutionary variable to be slightly non-solenoidal
while simultaneously setting up a mechanism to
drive the divergence residue towards zero.
More precisely,
we perturb the divergence-free evolutionary variable $\bm{u}$
in the time-derivative term 
to another velocity $\bm{w} := \bm{u} - \nabla\psi$
where $\psi$ is some scalar function;
meanwhile we change $\bm{u}$ to $\bm{w}$ in the diffusion term to
set up the mechanism of divergence decay.
Then there is no need to worry about
the influence of $\nabla\cdot\bm{w}\neq 0$
on numerical stability
because the evolution of $\bm{w}$ is not subject to
the divergence-free constraint.
These ideas lead to the GePUP formulation
proposed in \cite{zhang16:_GePUP},
which we briefly review here for completeness.

\subsection{The generic projection $\cProj$}
\label{sec:gener-proj-cproj}

A \emph{generic projection}
is a linear operator $\cProj$ on a vector space
satisfying
\begin{equation}
  \label{eq:cProj}
  \cProj \bm{v}^* = \bm{v} := \bm{v}^* -\nabla\psi,
\end{equation}
where $\psi$ is a scalar function
and $\nabla\cdot\bm{v}=0$ may or may not hold.
Since $\psi$ is not specified in terms of $\bm{v}^*$,
(\ref{eq:cProj}) is not a precise definition of $\cProj$, 
but rather a characterization of a family of operators,
which, in particular, includes
the Leray-Helmholtz projection $\cProjLH$. 
$\cProj$ can be used to perturb
the divergence-free velocity field $\bm{u}$
to some non-solenoidal velocity $\bm{w} := \cProj\bm{u}$
and is thus more flexible than $\cProjLH$
in characterizing discrete projections
that fail to fulfill the divergence-free constraint exactly.

It follows from (\ref{eq:cProj}) and
the commutativity of $\Delta$ and $\nabla$ that
\begin{equation}
  \label{eq:LaplProj}
  \Delta\cProj = \Delta - \nabla\nabla\cdot + \nabla\nabla\cdot\cProj.
\end{equation}

\subsection{GePUP: Reformulating INSE
  via generic projection and unconstrained PPE}
\label{sec:gepup:-reform-inse}

To accommodate the fact that
the discrete velocity field might not be divergence-free,
the evolutionary variable of the GePUP formulation is designed to be
a non-solenoidal velocity
\begin{equation}
  \label{eq:projectedVelocity}
  \bm{w} := \cProj\bm{u}
\end{equation}
instead of the divergence-free velocity $\bm{u}$
in the UPPE formulation (\ref{eq:UPPEstrong}).
The evolution of $\bm{w}$ is governed by
\begin{align*}
  \frac{\partial \bm{w}}{\partial t}
  &= \frac{\partial \cProj\bm{u}}{\partial t}
    = \cProj \frac{\partial \bm{u}}{\partial t}
    = \cProj \frac{\partial \bm{u}}{\partial t}
    - \nu(\Delta\bm{u} - \nabla\nabla\cdot\bm{u}
    + \nabla\nabla\cdot\cProj\bm{u} - \Delta\cProj\bm{u}) \\
  &= \cProj\bm{a} - \bm{a}^*
    + \bm{f} - \bm{u}\cdot\nabla\bm{u}
    + \nu\nabla\nabla\cdot\bm{u}
    - \nu\nabla\nabla\cdot\bm{w} + \nu\Delta\bm{w} \\
  &= \bm{f} - \bm{u}\cdot\nabla\bm{u} - \nabla q + \nu\Delta\bm{w},
\end{align*}
where the first step follows from (\ref{eq:projectedVelocity}),
the second from the commutativity of $\partial_t$ and $\cProj$,
the third from (\ref{eq:LaplProj}),
the fourth from (\ref{eq:projectedVelocity}) and
the definition of the Eulerian accelerations
in (\ref{eq:EulerianAccelerations}),
and the last from the definition
\begin{equation}
  \label{eq:gradQ}
  \nabla q := \bm{a}^* - \cProj\bm{a} - \nu\nabla\nabla\cdot\bm{u}
  +\nu\nabla\nabla\cdot\bm{w}.
\end{equation}
The RHS of (\ref{eq:gradQ}) is indeed a gradient field because of
the decomposition in (\ref{eq:INSEa-EulerianAccel}) and
the definition in (\ref{eq:cProj}).

Although the projected velocity $\bm{w}$ is allowed to be non-solenoidal, 
 we demand 
\begin{enumerate}[(i)]
\item
  $\bm{w} = \bm{g}$ on $\partial\Omega$,
\item
  its divergence $\nabla\cdot\bm{w}$ is governed by a heat equation
  with no-flux boundary conditions:
  \begin{subequations}
    \label{eq:assumptionsOfDivOfProjVel}
    \begin{alignat}{2}
      \frac{\partial (\nabla\cdot\bm{w})}{\partial t}
      &= \nu\Delta(\nabla\cdot\bm{w}) &\quad& \text{in } \Omega, \\
      \bm{n}\cdot\nabla\nabla\cdot\bm{w} &= 0 && \text{on }
      \partial\Omega.
    \end{alignat}
  \end{subequations}
\end{enumerate}
Consequently, we obtain 
\emph{the GePUP formulation of the INSE with Dirichlet boundary conditions}
\begin{subequations}
  \label{eq:GePUP}  
  \begin{alignat}{2}
    \frac{\partial \bm{w}}{\partial t} - \nu\Delta \bm{w}
    &= \bm{f} - \bm{u}\cdot\nabla\bm{u} - \nabla q &\quad&
    \text{in } \Omega, \\
    \bm{w} &= \bm{g} && \text{on } \partial\Omega, \\
    \bm{u} &= \cProjLH\bm{w} && \text{in } \Omega, \\
    \bm{u}\cdot\bm{n} &= \bm{g}\cdot\bm{n} && \text{on } \partial\Omega, \\
    \Delta q &= \nabla\cdot(\bm{f} - \bm{u}\cdot\nabla\bm{u}) &&
    \text{in } \Omega, \\
    \bm{n}\cdot\nabla q
    &= \bm{n}\cdot \left(\bm{f} - \bm{u}\cdot\nabla\bm{u}
      - \nu\nabla\times\nabla\times\bm{u}
      - \frac{\partial \bm{g}}{\partial t}\right) && \text{on } \partial\Omega,
  \end{alignat}
\end{subequations}
where (\ref{eq:GePUP}e) follows from the divergence of (\ref{eq:gradQ}),
(\ref{eq:EulerianAccelerations}),
the commutativity of $\nabla\cdot$ and $\Delta$,
(\ref{eq:projectedVelocity}),
and the assumption (\ref{eq:assumptionsOfDivOfProjVel}a),
(\ref{eq:GePUP}f) from
the normal component of (\ref{eq:gradQ}) on $\partial\Omega$,
(\ref{eq:EulerianAccelerations}),
the assumption (\ref{eq:assumptionsOfDivOfProjVel}b),
and the identity
\begin{equation}
  \label{eq:curlCurl}
  \Delta - \nabla\nabla\cdot = -\nabla\times\nabla\times.
\end{equation}
Since the generic projection $\cProj$ only perturbs
the divergence-free velocity $\bm{u}$ by a gradient field
and any gradient field is in
the null space of the Leray-Helmholtz projection $\cProjLH$,
(\ref{eq:GePUP}c) clearly holds.

Since (\ref{eq:assumptionsOfDivOfProjVel}a) is a starting point
of our deriving GePUP from INSE,
we expect that it can be recovered from (\ref{eq:GePUP}).
Indeed, 
 the divergence of (\ref{eq:GePUP}a) and (\ref{eq:GePUP}e) yield
 (\ref{eq:assumptionsOfDivOfProjVel}a), 
 which, together with the boundary condition (\ref{eq:assumptionsOfDivOfProjVel}b),
 drives $\nabla\cdot\bm{w}$ towards zero.
Hence the GePUP formulation preserves the advantage of the UPPE formulation (\ref{eq:UPPEstrong})
that the velocity divergence decays exponentially.

Starting from the INSE,
 we obtain (\ref{eq:GePUP}) as the consequence
 of (\ref{eq:assumptionsOfDivOfProjVel})
 that dictates the exponential decay of $\nabla\cdot\mathbf{w}$.
On the other hand,
 the INSE can be recovered from (\ref{eq:GePUP}) 
 by imposing the initial condition
\begin{equation}
  \label{eq:wInitialCondition}
  \bm{w}(t_0) = \bm{u}(t_0) \qquad \text{ in } \overline{\Omega},
\end{equation}
where $\overline{\Omega}$ is the closure of $\Omega$.
(\ref{eq:wInitialCondition}) is natural because, 
 designed as a perturbed version of $\bm{u}$,
 the non-solenoidal velocity $\bm{w}$ should converge to $\bm{u}$.
Therefore,
 \emph{the GePUP formulation is equivalent to the INSE}.

We summarize prominent features of the GePUP formulation (\ref{eq:GePUP})
as follows.
\begin{enumerate}[(GPF-1)]
\item
  The sole evolutionary variable
  is the non-solenoidal velocity $\bm{w}$,
  with $\bm{u}$ determined from $\bm{w}$
  via (\ref{eq:GePUP}c,d)
  and $q$ from $\bm{u}$
  via (\ref{eq:GePUP}e,f).
  This chain of determination
  $\bm{w}\rightarrow \bm{u} \rightarrow q$
  from Neumann BVPs
  is \emph{instantaneous}
  and has nothing to do with time integration.
  Therefore,
  a time integrator in MOL can be employed
  in a black-box manner
  and thus may be easily changed without affecting
  other modules of the entire solver.
\item
  There is no ambiguity on which velocities
  should be projected and which should not in MOL.
\item Now that the main evolutionary variable
  is formally non-solenoidal,
  the Leray-Helmholtz projection $\cProjLH$
  only comes into the system (\ref{eq:GePUP}) on the RHS.
  Although still contributing to the local truncation error,
  the approximation error of a discrete projection
  to $\cProjLH$ does not affect numerical stability. 
\end{enumerate}
The two issues with the UPPE formulation (\ref{eq:UPPEstrong})
discussed in the last two paragraphs
of Section \ref{sec:form-unconstr-ppe}
are resolved by (GPF-2) and (GPF-3).

\section{GePUP-FEM}
\label{sec:numerical-algorithms}

We develop high-order INSE solvers based on MOL,
where the governing equations are first discretized in space
by finite element methods
and then a high-order time integrator is applied to
the resulting system of ODEs.
To this end,
we first derive the weak form of the GePUP formulation (\ref{eq:GePUP}).

\subsection{The weak form of the GePUP formulation}
\label{sec:weak-formulation}

Denote by $\bm{H}_{\bm{g}}^1(\Omega) :=
\left\{ \bm{v}\in \bm{H}^1(\Omega):
  \bm{v} = \bm{g} \text{ on } \partial\Omega \right\}$,
multiply the GePUP formulation (\ref{eq:GePUP}) by test functions,
integrate by parts,
and we obtain the following \emph{weak form of
  the GePUP formulation (\ref{eq:GePUP})}.
\begin{enumerate}[(i)]
\item
  advance the velocity $\bm{w}$ by the momentum equation:
  find $\bm{w}\in L^2(0, T; \bm{H}^1_{\bm{g}}(\Omega))$
  such that
  \begin{equation}
    \label{eq:wWeakForm}
    \forall\bm{v}\in \bm{H}^1_{\bm{0}}(\Omega), \quad
    \left( \frac{\partial \bm{w}}{\partial t}, \bm{v} \right)_{\Omega}
    + \nu(\nabla\bm{w}, \nabla\bm{v})_{\Omega}
    = (\bm{f} - \bm{u}\cdot\nabla\bm{u} - \nabla q, \bm{v})_{\Omega};
  \end{equation}
\item
  compute the divergence-free velocity $\bm{u}$:
  find $\phi\in H^1_{\int=0}(\Omega)$ such that
  \begin{equation}
    \label{eq:uWeakForm}
    \forall\eta\in H^1_{\int=0}(\Omega), \quad
    (\nabla\phi, \nabla\eta)_{\Omega}
    = (\bm{w}, \nabla\eta)_{\Omega}
    - (\bm{n}\cdot\bm{g}, \eta)_{\partial\Omega},
  \end{equation}
  and set $\bm{u} = \bm{w} - \nabla\phi$;
\item
  extract the pressure $q$ from the velocity $\bm{u}$:
  find $q\in H^1_{\int=0}(\Omega)$ such that
  \begin{equation}
    \label{eq:qWeakForm}
    \forall\eta\in H^1_{\int=0}(\Omega), \quad
    (\nabla q, \nabla \eta)_{\Omega}
    = (\bm{f} - \bm{u}\cdot\nabla\bm{u}, \nabla\eta)_{\Omega}
    + \nu (\nabla\times\bm{u}, \bm{n}\times\nabla\eta)_{\partial\Omega}
    - \left(\bm{n}\cdot\frac{\partial\bm{g}}{\partial t}, \eta\right)_{\partial\Omega}.
  \end{equation}
\end{enumerate}
The second term on the RHS of (\ref{eq:qWeakForm}) is derived as follows.
\begin{align*}
  & (\bm{n}\cdot(\nabla\times\nabla\times\bm{u}), \eta)_{\partial\Omega}
    = (\nabla\times\nabla\times\bm{u}, \nabla\eta)_{\Omega}
    = (\nabla\times\nabla\times\bm{u}, \nabla\eta)_{\Omega}
    - (\nabla\times\bm{u}, \nabla\times\nabla\eta)_{\Omega} \\
  =\,& \int_{\Omega}\nabla\cdot \left( (\nabla\times\bm{u}) \times\nabla\eta \right)
     = \int_{\partial\Omega}\bm{n}\cdot \left( (\nabla\times\bm{u}) \times\nabla\eta \right)
     = -\left( \nabla\times\bm{u}, \bm{n}\times\nabla\eta \right)_{\partial\Omega},
\end{align*}
where the first step follows from the divergence theorem and
the identity (\ref{eq:curlCurl}),
the second from the fact that
the curl of any gradient field equals zero,
the third from the vector identity
\begin{equation*}
  \nabla\cdot(\bm{F}\times\bm{G})
  = (\nabla\times\bm{F})\cdot\bm{G}
  - (\nabla\times\bm{G})\cdot\bm{F},
\end{equation*}
the fourth from the divergence theorem,
and the last from the vector identity
\begin{equation*}
  \bm{n}\cdot (\bm{F}\times\bm{G}) = -\bm{F}\cdot (\bm{n}\times\bm{G}).
\end{equation*}

\revise{
Note that since the pressure is extracted by the weak form of an elliptic BVP,
we require it to have the $H^1$ regularity.
This differs from mixed formulations of the INSE,
where the $L^2$ regularity suffices for the pressure.
}

\subsection{Finite-element spatial discretization}
\label{sec:finite-elem-spat}

Let $\mathcal{T}_h$ be a mesh of the computational domain $\Omega$
with $h = \max_{K\in \mathcal{T}_h}h_K$,
$\bm{V}_h\subset \bm{H}^1(\Omega)$
a space of $C^0$ finite elements
for approximating the velocity field,
$\bm{V}_{0, h} = \bm{V}_h\cap \bm{H}^1_0(\Omega)$
the subspace of $\bm{V}_h$
consisting of vector fields that vanish on $\partial\Omega$,
and $Q_h\subset H^1_{\int=0}(\Omega)$
a space of $C^0$ finite elements
for approximating the pressure field.
For domains with curved boundaries,
isoparametric finite elements are employed.

Because the pressure is explicitly determined by
the Poisson equation (\ref{eq:qWeakForm}),
continuous Lagrange finite elements of equal order can be used
for $\bm{V}_h$ and $Q_h$.
This contrasts with mixed finite element methods,
where employing equal-order elements leads to instabilities
due to the violation of the inf-sup condition (\ref{eq:infSupCondition}).

Let $\left\{ \eta_j \right\}_{j=1}^N$ be
a basis of the continuous Lagrange finite element space
and express the solutions as
\begin{equation}
  \label{eq:solutionFEExpansions}
  \begin{aligned}
    \bm{w}_h(\bm{x}, t) &= \sum\nolimits_{d=1}^{\Dim}w_{h, d}(\bm{x}, t)\bm{e}_d,
    &\ &
    w_{h, d}(\bm{x}, t) = \sum\nolimits_{j=1}^Nw_{d, j}(t)\eta_j(\bm{x}),
    &\ &
    \phi_h(\bm{x}, t) = \sum\nolimits_{j=1}^N \phi_j(t)\eta_j(\bm{x}),
    \\
    \bm{u}_h(\bm{x}, t) &= \sum\nolimits_{d=1}^{\Dim}u_{h, d}(\bm{x}, t)\bm{e}_d,
    &\ &
    u_{h, d}(\bm{x}, t) = \sum\nolimits_{j=1}^Nu_{d, j}(t)\eta_j(\bm{x}),
    &&
    q_h(\bm{x}, t) = \sum\nolimits_{j=1}^Nq_j(t)\eta_j(\bm{x}),
  \end{aligned}
\end{equation}
where the subscript $_d$ denotes a vector component,
$\{\bm{e}_d\}_{d=1}^{\Dim}$ is
the standard basis of $\mathbb{R}^{\Dim}$,
$w_{d,j}(t), \phi_j(t), u_{d,j}(t)$, and $q_j(t)$ are
the unknown coefficients to be sought.
Restricting the weak formulations
(\ref{eq:wWeakForm})--(\ref{eq:qWeakForm})
to the finite-dimensional spaces
$\bm{V}_{h, 0}$ and $Q_h$ and
substituting (\ref{eq:solutionFEExpansions}) into the resulting equations,
we obtain the following semi-discrete algorithm:
\begin{subequations}
  \label{eq:spatialDiscretization}
  \begin{align}
    M\frac{\dif \bm{W}_{h, d}(t)}{\dif t} + \nu A\bm{W}_{h, d}(t)
    &= \bm{F}^{\bm{w}}_d(\bm{u}_h, q_h, \bm{f}), \\
    A\bm{\Phi}_h &= \bm{F}^{\phi}(\bm{w}_h, \bm{g}), \\
    M\bm{U}_{h, d} &= \bm{F}^{\bm{u}}_d(w_{h, d}, \phi_h), \\
    A\bm{Q}_h &= \bm{F}^q(\bm{u}_h, \bm{f}, \bm{g}),
  \end{align}
\end{subequations}
where
$\bm{W}_{h, d}\in \mathbb{R}^N$,
$\bm{\Phi}_h\in \mathbb{R}^N$,
$\bm{U}_{h, d}\in \mathbb{R}^N$,
and $\bm{Q}_h\in \mathbb{R}^N$ are
the vectors with components
$w_{d, j}(t)$, $\phi_j(t)$, $u_{d, j}(t)$, and $q_j(t)$,
respectively,
$M\in \mathbb{R}^{N\times N}$
and $A\in \mathbb{R}^{N\times N}$ are
the mass and stiffness matrices,
i.e.,
$m_{ij} = (\eta_i, \eta_j)_{\Omega}$ and
$a_{ij} = (\nabla\eta_i, \nabla\eta_j)_{\Omega}$,
and the $i$th components of the RHS of (\ref{eq:spatialDiscretization}) are
respectively
\begin{subequations}
  \label{eq:spatialDiscretizationRHS}
  \begin{align}
    F^{\bm{w}}_{d, i}(\bm{u}_h, q_h, \bm{f})
    &= \left( f_d - \bm{u}_h\cdot\nabla u_{h, d} - \frac{\partial q_h}{\partial x_d}, \eta_i \right)_{\Omega}, \\
    F^{\phi}_i(\bm{w}_h, \bm{g})
    &= \left( \bm{w}_h, \nabla\eta_i \right)_{\Omega} - \left( \bm{n}\cdot\bm{g}, \eta_i \right)_{\partial\Omega}, \\
    F^{\bm{u}}_{d, i}(w_{h, d}, \phi_h)
    &= \left( w_{h, d} - \frac{\partial\phi_h}{\partial x_d}, \eta_i \right)_{\Omega}, \\
    F^q_i(\bm{u}_h, \bm{f}, \bm{g})
    &= \left( \bm{f} - \bm{u}_h\cdot\nabla\bm{u}_h, \nabla\eta_i \right)_{\Omega}
      + \nu \left( \nabla\times\bm{u}_h, \bm{n}\times\nabla\eta_i \right)_{\partial\Omega}
      - \left( \bm{n}\cdot \frac{\partial \bm{g}}{\partial t}, \eta_i \right)_{\partial\Omega},
  \end{align}
\end{subequations}
\revise{
where the temporal derivative term $\frac{\partial \bm{g}}{\partial t}$
in (\ref{eq:spatialDiscretizationRHS}d)
can be computed using the analytic expression of
the explicitly provided velocity Dirichlet boundary function $\bm{g}(\bm{x}, t)$.
}

\subsection{Temporal integration}
\label{sec:temporal-integration}

Besides the explicit treatment of the pressure gradient term to
avoid solving large coupled saddle-point systems,
it is also desirable
for flows with low to moderate Reynolds numbers
to treat the convection term explicitly and
the diffusion term implicitly,
because implicit treatment of the non-stiff convection term leads to
a nonlinear system challenging for iterative solvers
whereas explicit treatment of the stiff diffusion term
incurs a stringent time step constraint.
For this purpose,
an implicit-explicit Runge-Kutta scheme
\cite{ascher97:_implic_runge_kutta}
can be seamlessly integrated with the semi-discrete GePUP algorithm
(\ref{eq:spatialDiscretization})
to form a group of semi-implicit time-marching schemes for the INSE.
In this work,
we adopt two ERK-ESDIRK methods:
one is the six-stage, fourth-order accurate
ARK4(3)6L[2]SA scheme
proposed in \cite{kennedy03:_addit_runge_kutta_schem_for};
the other is the eight-stage, fifth-order accurate
ARK5(4)8L[2]$\text{SA}_2$ scheme
proposed in \cite{kennedy19:_higher_runge_kutta}.
Both are L-stable and stiffly accurate.

For an ODE of the form
\begin{displaymath}
  \frac{\dif \zeta}{\dif t} = \exOp(\zeta, t) + \imOp(\zeta),
\end{displaymath}
the steps of an ERK-ESDIRK method are
\begin{subequations}
  \label{eq:ark4modelODE}
  \begin{align}
    \zeta^{(1)} &= \zeta^n \approx \zeta(t^n),
    \\
    \forall s = 2, 3, \ldots, \nStages, \quad
    \left(I - \dt\gamma\imOp \right) \zeta^{(s)}
                &= \zeta^{n}
                  + \dt \sum\limits^{s-1}_{j=1}  a^{[\text{E}]}_{s, j} \exOp
                  \left(\zeta^{(j)}, t^{(j)}\right)
                  + \dt \sum\limits^{s-1}_{j=1}  a^{[\text{I}]}_{s, j} \imOp\left(\zeta^{(j)}\right),
    \\
    \zeta^{n+1}
                &= \zeta^{(\nStages)}
                  + \dt \sum^{\nStages}_{j=1} \left(b_j - a^{[\text{E}]}_{\nStages, j}\right)
                  \exOp
                  \left(\zeta^{(j)}, t^{(j)}\right),
  \end{align}
\end{subequations}
where the superscript $^{(s)}$ denotes an intermediate stage,
$t^{(s)}=t^n+c_s \dt$ the time of that stage,
$\nStages$ the number of stages,
$A^{[\text{E}]}, \bm{b}$ and $\bm{c}$
the coefficients of the Butcher tableau of
the underlying ERK method,
and $A^{[\text{I}]}, \bm{b}$ and $\bm{c}$
those of the underlying ESDIRK method.
Note that (\ref{eq:ark4modelODE}c) follows from
\begin{equation*}
  \zeta^{n+1} = \zeta^n
  + \dt\sum_{j=1}^{\nStages}b_j^{[\text{E}]}\exOp\left( \zeta^{(j)}, t^{(j)} \right)
  + \dt\sum_{j=1}^{\nStages}b_j^{[\text{I}]}\imOp\left( \zeta^{(j)} \right)
\end{equation*}
and $b^{[\text{E}]}_{j}=b^{[\text{I}]}_{j}=b_j=a^{[\text{I}]}_{\nStages,j}$.
For more details on this method,
the reader is referred to
\cite{kennedy03:_addit_runge_kutta_schem_for,
  kennedy19:_higher_runge_kutta}.

Applying the ERK-ESDIRK method (\ref{eq:ark4modelODE}) to
the ODE system (\ref{eq:spatialDiscretization}) with
the splitting
\begin{equation*}
  \exOp
  := \bm{F}_d^{\bm{w}}(\bm{u}_h, q_h, \bm{f}),
  \qquad
  \imOp
  := -\nu A\bm{W}_{h, d}
\end{equation*}
yields the following fully discrete GePUP-FEM algorithm:
\begin{subequations}
  \label{eq:GePUP-IMEX}
  \begin{align}
    & \qquad \qquad \qquad \qquad \qquad \qquad \qquad \quad\,\,\,
      \bm{W}_{h, d}^{(1)} = \bm{W}_{h, d}^n,
    \\
    &
      \forall s = 2, 3, \ldots, \nStages, \quad
      \left\{
      \begin{aligned}
        \bigl(M + \nu\dt\gamma A\bigr) \bm{W}_{h, d}^{(s)}
        &=
        M\bm{W}_{h, d}^n
        +\dt \sum\limits^{s-1}_{j=1}  a^{[\text{E}]}_{s,j}
        \bm{F}^{\bm{w}}_d \left( \bm{u}_h^{(j)}, q_h^{(j)}, \bm{f}^{(j)} \right)
        + \nu\dt \sum\limits^{s-1}_{j=1}  a^{[\text{I}]}_{s,j} A\bm{W}_{h, d}^{(j)},
        \\
        A\bm{\Phi}_h^{(s)}
        &= \bm{F}^{\phi}\left( \bm{w}_h^{(s)}, \bm{g}^{(s)} \right),
        \\
        M\bm{U}_{h, d}^{(s)}
        &= \bm{F}^{\bm{u}}_d \left( w_{h, d}^{(s)}, \phi_h^{(s)} \right),
        \\
        A\bm{Q}_h^{(s)}
        &= \bm{F}^q \left( \bm{u}_h^{(s)}, \bm{f}^{(s)}, \bm{g}^{(s)} \right),
      \end{aligned}
          \right.
    \\
    & \qquad \qquad \qquad \qquad \qquad \qquad  \quad\,\,
    \left\{
      \begin{aligned}
        M\bm{W}_{h, d}^{*} &= M\bm{W}_{h, d}^{(\nStages)}
        + \dt \sum\limits^{\nStages}_{j=1} \left(b_j - a^{[\text{E}]}_{\nStages,j}\right)
        \bm{F}^{\bm{w}}_d \left(\bm{u}_h^{(j)}, q_h^{(j)}, \bm{f}^{(j)}\right),
        \\
        A\bm{\Phi}_h^{n+1} &=
        \bm{F}^{\phi}\left( \bm{w}_h^{*}, \bm{g}^{n+1} \right),
        \\
        M\bm{U}_{h, d}^{n+1}
        &= \bm{F}^{\bm{u}}_d \left( w_{h, d}^{*}, \phi_h^{n+1} \right),
        \\
        A\bm{Q}_h^{n+1}
        &= \bm{F}^q \left( \bm{u}_h^{n+1}, \bm{f}^{n+1}, \bm{g}^{n+1} \right),
        \\
        \bm{W}_{h, d}^{n+1} &= \bm{U}_{h, d}^{n+1}.
      \end{aligned}
    \right.
  \end{align}
\end{subequations}
Due to the nonlinear convection term in $\exOp$,
the divergence of $\bm{w}_h^{*}$ in (\ref{eq:GePUP-IMEX}c)
is much greater than
that of $\bm{w}_h^{(s)}$ in previous stages.
Hence,
the final solution $\bm{w}_h^{n+1}$ of this time step is set to
$\bm{u}_h^{n+1}$,
not $\bm{w}_h^*$.

In summary,
our numerical algorithms are
\begin{itemize}
\item
  versatile:
  time integrators can be employed in a black-box manner;
\item
  accurate:
  fourth- and higher-order accuracy both in time and in space is achieved;
\item
  efficient:
  at each time step,
  advancing the INSE amounts to solving a sequence of elliptic BVPs,
  for which highly efficient linear solvers are readily available;
\item
  flexible:
  the selected finite element spaces for the velocity and the pressure 
  do not have to satisfy the inf-sup condition (\ref{eq:infSupCondition}).
\end{itemize}

\section{Numerical tests}
\label{sec:numerical-tests}

In this section,
we test GePUP-FEM by standard benchmark problems to
demonstrate its fourth- and higher-order accuracy and
its capability of accurately and efficiently resolving the right physics.
Our numerical algorithm is implemented
using the state-of-the-art open-source \texttt{C++} finite element library
\textsf{deal.II}
\cite{
  arndt21},
which supports
distributed quadrilateral/hexahedral meshes
via \textsf{p4est} \cite{burstedde11},
matrix-free computations with sum factorization
\cite{kronbichler12},
and adaptive mesh refinement (AMR)
with hanging nodes
\cite{bangerth11:_algor}.
The resulting linear systems are solved by
the conjugate gradient method
preconditioned by geometric multigrid
\cite{clevenger21:_flexib_paral_adapt_geomet_multig_method_fem}.
All our numerical simulations were performed on a server
with AMD Ryzen Threadripper PRO 3995WX 64-Core CPUs
running at 2.70 GHz.

As discussed in the previous section,
we employ continuous Lagrange finite elements of equal order
for the velocity and pressure finite element spaces
which are known to violate the inf-sup condition
(\ref{eq:infSupCondition}).
The time step size $\dt$ is determined by the Courant number
\begin{equation}
  \label{eq:courantNumber}
  \text{Cr} := k\dt \max_{K\in \mathcal{T}_h}
  \frac{\|\bm{u}_h\|_{\bm{L}^{\infty}(K)}}{h_K},
\end{equation}
where $k$ is the polynomial degree
of the employed Lagrange finite element.

\subsection{Numerical tests for convergence rates}
\label{sec:convergence-tests}

In this subsection,
standard benchmark problems are solved to
verify the high-order accuracy of GePUP-FEM both in time and in space.
We test two combinations:
one is a fourth-order GePUP-FEM scheme
with ARK4(3)6L[2]SA as the time integrator and
continuous $\mathbb{Q}_3$-elements for
both the velocity and pressure finite element spaces;
the other is a fifth-order GePUP-FEM scheme
with ARK5(4)8L[2]$\text{SA}_2$ as the time integrator and
continuous $\mathbb{Q}_4$-elements for
both the velocity and pressure finite element spaces.

\subsubsection{Taylor-Green vortex in two dimensions}
\label{sec:taylor-green-vortex}

This test concerns the two-dimensional Taylor-Green vortex
\cite{taylor37:_mechan,
  nguyen11:_galer_navier_stokes}
with analytic solutions
\begin{subequations}
  \label{eq:taylorGreenVortex}
  \begin{align}
    \bm{u}(x, y, t) &= e^{-\frac{2\pi^2 t}{\text{Re}}}
                      \begin{pmatrix}
                        -\cos(\pi x)\sin(\pi y) \\ \sin(\pi x)\cos(\pi y)
                      \end{pmatrix}, \\
    p(x, y, t) &= -\frac{1}{4}e^{-\frac{4\pi^2t}{\text{Re}}}
                 \left( \cos(2\pi x) + \cos(2\pi y) \right),
  \end{align}
\end{subequations}
where $\text{Re} = 1/\nu$ is the Reynolds number.
The time derivative of the velocity cancels the diffusion term and
the pressure gradient cancels the convection term,
resulting in a zero external force.
The Dirichlet boundary condition for the velocity is taken as
the restriction of the exact solution to the domain boundary and
the initial condition as an instantiation of the exact solution at $t_0$.

In the unit square $\Omega = (0, 1)^2$,
(\ref{eq:taylorGreenVortex}) is advanced
from $t_0 = 0$ to $t_{\mathrm{e}} = 1$ with Cr = 0.8
on four successively refined square grids with size $h$.
Errors and convergence rates of
the fourth- and fifth-order GePUP-FEM schemes
for $\text{Re} = 100$
are summarized in Tables
\ref{tab:taylorGreenVortexRe1e2FourthOrder} and
\ref{tab:taylorGreenVortexRe1e2FifthOrder},
respectively,
which clearly show the expected convergence rates
for both the velocity and the pressure,
i.e.,
for the $k$th-order GePUP-FEM scheme,
both the velocity and the pressure are $k$th-order accurate
in the $L^2$- and $L^{\infty}$-norms and
$(k-1)$th-order accurate in the $H^1$-norm.

\begin{table}
  \caption{Errors and convergence rates of
    the fourth-order GePUP-FEM scheme for the Taylor-Green vortex
    with Re = 100,
    $t_0 = 0.0$,
    $t_{\mathrm{e}} = 1.0$,
    and Cr = 0.8.}
  \centering
  \renewcommand{\arraystretch}{1.3}
  \begin{tabular}{ccccccccc}
  \hline
  \multicolumn{2}{c}{} & $h=\frac{1}{8}$ & Rate & $h=\frac{1}{16}$ & Rate & $h=\frac{1}{32}$ & Rate & $h=\frac{1}{64}$
  \\ \hline \hline
  \multirow{3}*{$\bm{u}$} & $\bm{L}^2$ & 9.28e-06 & 4.65 & 3.69e-07 & 4.43 & 1.71e-08 & 4.19 & 9.38e-10
  \\ & $\bm{H}^1$ & 8.47e-04 & 3.44 & 7.80e-05 & 3.15 & 8.77e-06 & 3.04 & 1.07e-06
  \\ & $\bm{L}^{\infty}$ & 4.19e-05 & 4.47 & 1.89e-06 & 4.24 & 9.97e-08 & 4.17 & 5.53e-09
  \\ \hline
  
  \multirow{3}*{$q$} & $L^2$ & 1.68e-05 & 4.03 & 1.03e-06 & 4.02 & 6.36e-08 & 4.01 & 3.95e-09
  \\ & $H^1$ & 1.64e-03 & 3.01 & 2.03e-04 & 3.00 & 2.54e-05 & 3.00 & 3.17e-06
  \\ & $L^{\infty}$ & 5.57e-05 & 3.94 & 3.62e-06 & 3.99 & 2.29e-07 & 3.97 & 1.46e-08
  \\ \hline
\end{tabular}

  \label{tab:taylorGreenVortexRe1e2FourthOrder}
\end{table}

\begin{table}
  \caption{Errors and convergence rates of
    the fifth-order GePUP-FEM scheme for the Taylor-Green vortex
    with Re = 100,
    $t_0 = 0.0$,
    $t_{\mathrm{e}} = 1.0$,
    and Cr = 0.8.}
  \centering
  \renewcommand{\arraystretch}{1.3}
  \begin{tabular}{ccccccccc}
  \hline
  \multicolumn{2}{c}{} & $h=\frac{1}{8}$ & Rate & $h=\frac{1}{16}$ & Rate & $h=\frac{1}{32}$ & Rate & $h=\frac{1}{64}$
  \\ \hline \hline
  \multirow{3}*{$\bm{u}$} & $\bm{L}^2$ & 3.90e-07 & 5.73 & 7.35e-09 & 5.67 & 1.44e-10 & 5.47 & 3.26e-12
  \\ & $\bm{H}^1$ & 5.66e-05 & 4.70 & 2.18e-06 & 4.69 & 8.46e-08 & 4.50 & 3.74e-09
  \\ & $\bm{L}^{\infty}$ & 1.24e-06 & 5.61 & 2.53e-08 & 5.69 & 4.88e-10 & 5.55 & 1.04e-11
  \\ \hline
  
  \multirow{3}*{$q$} & $L^2$ & 6.56e-07 & 4.97 & 2.09e-08 & 4.83 & 7.37e-10 & 4.54 & 3.16e-11
  \\ & $H^1$ & 8.29e-05 & 4.03 & 5.06e-06 & 4.01 & 3.13e-07 & 4.00 & 1.95e-08
  \\ & $L^{\infty}$ & 2.34e-06 & 4.92 & 7.71e-08 & 4.97 & 2.47e-09 & 4.79 & 8.91e-11
  \\ \hline
\end{tabular}

  \label{tab:taylorGreenVortexRe1e2FifthOrder}
\end{table}

\subsubsection{Beltrami flow in three dimensions}
\label{sec:beltrami-flow}

In this test,
we consider the three-dimensional Beltrami flow
\cite{ethier94:_exact_navier_stokes,
  piatkowski18:_galer_navier_stokes}
with exact solutions
\begin{subequations}
  \label{eq:beltramiFlow}
  \begin{align}
    \bm{u}(x, y, z, t) =& -a e^{-\frac{d^2t}{\text{Re}}}
    \begin{pmatrix}
      e^{ax}\sin(ay + dz) + e^{az}\cos(ax + dy) \\
      e^{ay}\sin(az + dx) + e^{ax}\cos(ay + dz) \\
      e^{az}\sin(ax + dy) + e^{ay}\cos(az + dx)
    \end{pmatrix}, \\
    p(x, y, z, t) =& \nonumber -\frac{a^2}{2}e^{-\frac{2d^2t}{\text{Re}}}
    \Big( e^{2ax} + e^{2ay} + e^{2az} + 2\sin(ax + dy)\cos(az + dx)e^{a(y+z)} \\
    &+ 2\sin(ay + dz)\cos(ax + dy)e^{a(z + x)}
    + 2a\sin(az + dx)\cos(ay + dz)e^{a(x + y)} \Big),
  \end{align}
\end{subequations}
where $\text{Re} = 1/\nu$ is the Reynolds number,
$a = \pi/4$,
and $b = \pi/2$.
As in Section \ref{sec:taylor-green-vortex}, 
the time derivative of the velocity cancels the diffusion term and
the pressure gradient cancels the convection term,
resulting in a zero forcing term.
The Dirichlet boundary condition for the velocity is taken as
the restriction of the exact solution to the domain boundary and
the initial condition as an instantiation of the exact solution at $t_0$.

In the cube $\Omega = (-1, 1)^3$,
(\ref{eq:beltramiFlow}) is advanced from
$t_0 = 0$ to $t_{\mathrm{e}} = 1$ with Cr = 0.4
on four successively refined uniform hexahedral meshes
with element side length $h$.
Errors and convergence rates of
the fourth- and fifth-order GePUP-FEM schemes
for Re = 100
are presented in Tables
\ref{tab:beltramiFlowRe1e2FourthOrder} and
\ref{tab:beltramiFlowRe1e2FifthOrder},
respectively.
It is observed that
while the velocity exhibits the expected convergence rates in all norms,
the pressure has order reductions
in the $L^2$- and $L^{\infty}$-norms.
These order reductions may be caused by
the evaluation of the term
$(\nabla\times\bm{u}_h, \bm{n}\times\nabla\eta_i)_{\partial\Omega}$
in (\ref{eq:spatialDiscretizationRHS}d)
and the fact that the direct evaluation of $\nabla\times\bm{u}_h$
from the finite element solution $\bm{u}_h$ is of sub-optimal order of accuracy.
However,
we see that
the convergence rate of the pressure in the $H^1$-norm is as expected.

\revise{
To demonstrate the efficiency advantage of GePUP-FEM,
we compare it with a second-order consistent splitting (CS) scheme
\cite{guermond06:_overview_of_projection_methods_for_incompressible_flows}
by their running times in achieving the same accuracy.
For this comparison,
we implemented the CS scheme using
a second-order implicit-explicit Runge-Kutta method
\cite{ascher97:_implic_runge_kutta}
as the time integrator and
continuous $\mathbb{Q}_1$-elements
for both the velocity and pressure finite element spaces.
Errors and convergence rates of this second-order CS scheme
are shown in Table \ref{tab:beltramiFlowRe1e2SecondOrder}.
The CPU times (in seconds) for the second-order CS scheme are
$0.22$, $1.85$, $24.73$, and $307.53$
on the four successively refined meshes,
while those for the fourth-order GePUP-FEM scheme are
$5.04$, $15.92$, $101.16$, and $1258.68$.
The $\bm{L}^2$-norm of the velocity error
 of the fourth-order GePUP-FEM scheme 
 with $h = 2^{-2}$ 
 is $4.05\times 10^{-5}$
 while that of the second-order CS scheme
 with $h = 2^{-4}$
 is $5.10\times 10^{-3}$. 
For the second-order CS scheme
 to achieve the same accuracy $4.05\times 10^{-5}$
 of the fourth-order GePUP-FEM scheme,
 even with an optimal complexity of solving linear systems, 
 the mesh has to be refined
 at least 
 $\log_{2^2}\frac{5.10\times 10^{-3}}{4.05\times 10^{-5}} \approx 3.49$ times,
 resulting in a running time of
 at least $307.53 \times (2^3 \times 2)^{3.49} \approx 4.9\times 10^6$ seconds,
 where the ``$\times 2$'' comes from the fact that
 the time step size is also halved during mesh refinement.
Thus,
the speedup of the fourth-order GePUP-FEM scheme
over the second-order CS scheme
in achieving a velocity $\bm{L}^2$ accuracy of $4.05\times 10^{-5}$ is
$4.9\times 10^6 / 15.92 \approx 3.08\times 10^5$!
This clearly illustrates the efficiency advantage of
the high-order accurate GePUP-FEM schemes.
}

\begin{table}
  \caption{Errors and convergence rates of
    the fourth-order GePUP-FEM scheme for the Beltrami flow
    with Re = 100,
    $t_0 = 0.0$,
    $t_{\mathrm{e}} = 1.0$,
    and Cr = 0.4.}
  \centering
  \renewcommand{\arraystretch}{1.3}
  \begin{tabular}{ccccccccc}
  \hline
  \multicolumn{2}{c}{} & $h=\frac{1}{2}$ & Rate & $h=\frac{1}{4}$ & Rate & $h=\frac{1}{8}$ & Rate & $h=\frac{1}{16}$
  \\ \hline \hline
  \multirow{3}*{$\bm{u}$} & $\bm{L}^2$ & 8.91e-04 & 4.46 & 4.05e-05 & 4.12 & 2.32e-06 & 4.16 & 1.30e-07
  \\ & $\bm{H}^1$ & 1.76e-02 & 3.37 & 1.70e-03 & 3.12 & 1.96e-04 & 3.08 & 2.32e-05
  \\ & $\bm{L}^{\infty}$ & 1.12e-03 & 4.41 & 5.27e-05 & 4.12 & 3.03e-06 & 4.08 & 1.79e-07
  \\ \hline

  \multirow{3}*{$q$} & $L^2$ & 1.00e-03 & 4.14 & 5.68e-05 & 3.89 & 3.84e-06 & 3.61 & 3.16e-07
  \\ & $H^1$ & 1.99e-02 & 3.13 & 2.28e-03 & 3.02 & 2.81e-04 & 3.01 & 3.49e-05
  \\ & $L^{\infty}$ & 1.93e-03 & 3.29 & 1.98e-04 & 3.59 & 1.64e-05 & 3.62 & 1.33e-06
  \\ \hline
\end{tabular}

  \label{tab:beltramiFlowRe1e2FourthOrder}
\end{table}

\begin{table}
  \caption{Errors and convergence rates of
    the fifth-order GePUP-FEM scheme for the Beltrami flow
    with Re = 100,
    $t_0 = 0.0$,
    $t_{\mathrm{e}} = 1.0$,
    and Cr = 0.4.}
  \centering
  \renewcommand{\arraystretch}{1.3}
  \begin{tabular}{ccccccccc}
  \hline
  \multicolumn{2}{c}{} & $h=\frac{1}{2}$ & Rate & $h=\frac{1}{4}$ & Rate & $h=\frac{1}{8}$ & Rate & $h=\frac{1}{16}$
  \\ \hline \hline
  \multirow{3}*{$\bm{u}$} & $\bm{L}^2$ & 6.13e-05 & 5.15 & 1.73e-06 & 5.65 & 3.44e-08 & 5.63 & 6.93e-10
  \\ 
  & $\bm{H}^1$ & 1.80e-03 & 4.17 & 1.00e-04 & 4.59 & 4.16e-06 & 4.53 & 1.80e-07
  \\ 
  & $\bm{L}^{\infty}$ & 9.45e-05 & 4.84 & 3.30e-06 & 5.50 & 7.31e-08 & 4.66 & 2.88e-09
  \\ \hline

  \multirow{3}*{$q$} & $L^2$ & 7.70e-05 & 4.87 & 2.63e-06 & 5.02 & 8.15e-08 & 4.66 & 3.21e-09
  \\ 
  & $H^1$ & 2.41e-03 & 4.09 & 1.41e-04 & 4.39 & 6.74e-06 & 4.26 & 3.53e-07
  \\ 
  & $L^{\infty}$ & 2.03e-04 & 4.55 & 8.66e-06 & 4.02 & 5.36e-07 & 3.90 & 3.58e-08
  \\ \hline
\end{tabular}

  \label{tab:beltramiFlowRe1e2FifthOrder}
\end{table}

\begin{table}
  \caption{\revise{Errors and convergence rates of
    a second-order consistent splitting scheme
    \cite{guermond06:_overview_of_projection_methods_for_incompressible_flows}
    for the Beltrami flow
    with Re = 100,
    $t_0 = 0.0$,
    $t_{\mathrm{e}} = 1.0$,
    and Cr = 0.4.}}
  \centering
  \renewcommand{\arraystretch}{1.3}
  \revise{
    \begin{tabular}{ccccccccc}
  \hline
  \multicolumn{2}{c}{} & $h=\frac{1}{2}$ & Rate & $h=\frac{1}{4}$ & Rate & $h=\frac{1}{8}$ & Rate & $h=\frac{1}{16}$
  \\ \hline \hline
  \multirow{3}*{$\bm{u}$} & $\bm{L}^2$ & 3.33e-01 & 1.97 & 8.47e-02 & 2.03 & 2.07e-02 & 2.02 & 5.10e-03
  \\ & $\bm{H}^1$ & 1.98e+00 & 1.14 & 8.99e-01 & 1.06 & 4.33e-01 & 1.02 & 2.14e-01
  \\ & $\bm{L}^{\infty}$ & 2.34e-01 & 1.89 & 6.34e-02 & 2.07 & 1.51e-02 & 2.09 & 3.54e-03
  \\ \hline

  \multirow{3}*{$p$} & $L^2$ & 5.41e-01 & 1.87 & 1.47e-01 & 1.98 & 3.75e-02 & 1.97 & 9.59e-03
  \\ & $H^1$ & 2.04e+00 & 1.23 & 8.68e-01 & 1.10 & 4.06e-01 & 1.03 & 1.99e-01
  \\ & $L^{\infty}$ & 6.52e-01 & 1.10 & 3.05e-01 & 1.68 & 9.50e-02 & 1.78 & 2.77e-02
  \\ \hline
\end{tabular}

  }
  \label{tab:beltramiFlowRe1e2SecondOrder}
\end{table}

\subsection{Acceptance tests via physical benchmarks}
\label{sec:numerical-benchmarks}

In this subsection,
standard benchmark problems,
including a single-vortex test,
\revise{
a lid-driven cavity,
}
two-dimensional flows past a circular cylinder,
and three-dimensional flows past a sphere,
are simulated to demonstrate that
GePUP-FEM is capable of
accurately and efficiently resolving the right physics,
i.e.,
\emph{converging to the correct solution}.
We only test the fourth-order GePUP-FEM,
\revise{
i.e.,
ARK4(3)6L[2]SA as the time integrator and
continuous (isoparametric) $\mathbb{Q}_3$-elements
for both the velocity and pressure finite element spaces.
}

To efficiently resolve strongly localized features, 
 we augment GePUP-FEM to support adaptive mesh refinement
 via the steps 
$\texttt{SOLVE} \rightarrow \texttt{ESTIMATE} \rightarrow
\texttt{MARK} \rightarrow \texttt{REFINE}$.
In  the \texttt{ESTIMATE} step,
we calculate the indicator of local refinement by
\begin{equation*}
  \forall K\in \mathcal{T}_h, \quad
  \eta_K = h_K\|\nabla\times\bm{u}_h\|_{\bm{L}^{\infty}(K)}. 
\end{equation*}
In the \texttt{MARK} step,
 the D\"{o}rfler marking strategy
\cite{doerfler96:_poiss} is employed
to mark
the set $\mathcal{R}_h$ of all refined elements and
the set $\mathcal{C}_h$ of all coarsened elements as 
\begin{subequations}
  \label{eq:dorflerStrategy}
  \begin{align}
    \text{find } \mathcal{R}_h\subset \mathcal{T}_h
    \text{ s.t. }  \sum_{K\in\mathcal{R}_h}\eta_K \ge
    \theta_{\mathrm{R}} \sum_{K\in \mathcal{T}_h}\eta_K
    \text{ with minimal cardinality}, \\
    \text{find } \mathcal{C}_h\subset \mathcal{T}_h
    \text{ s.t. }  \sum_{K\in\mathcal{C}_h}\eta_K \le
    \theta_{\mathrm{C}} \sum_{K\in \mathcal{T}_h}\eta_K
    \text{ with maximal cardinality}, 
  \end{align}
\end{subequations}
where $\theta_{\mathrm{R}}\in (0, 1)$
and $\theta_{\mathrm{C}}\in (0, 1)$
are the refining threshold 
and the coarsening threshold, respectively.
The D\"orfler marking strategy (\ref{eq:dorflerStrategy})
is known \cite{pfeiler20:_doerf} to have optimal complexity.
For more details on the algorithms and data structures
for massively parallel adaptive finite element codes,
the reader is referred to \cite{bangerth11:_algor}.

\subsubsection{Single vortex}
\label{sec:single-vortex}

While the tests in Section \ref{sec:convergence-tests} focus on
the high-order accuracy of our numerical solver,
the purpose of this test is to
illustrate its efficiency in accurately capturing
\revise{
sharp gradients and localized features
near the boundaries,
which are characteristic of boundary layer-like structures in this 2D example.
}
Following \cite{bell91},
we consider an axisymmetric velocity field
defined on the unit square $(0, 1)^2$ by
\begin{equation*}
  u_{\theta}(r_v) =
  \begin{cases}
    \frac{1}{2}r_v - 4r_v^3 & \text{if } r_v < R; \\
    \frac{R}{r_v}\left( \frac{1}{2}R - 4R^3 \right) & \text{if } r_v \ge R,
  \end{cases}
\end{equation*}
where $r_v$ is the distance from $(x, y)$
to the vortex center $(0.5, 0.5)$ and
$R = 0.2$.
Before being used as the initial condition of the INSE,
this non-solenoidal vector field is projected by
the discrete Leray-Helmholtz projection operator
so that it is approximately divergence-free.
The velocity satisfies no-slip boundary conditions
$\bm{u} = \bm{0}$ on the whole domain boundary.
The external forcing function $\bm{f}$ is taken to be zero.
The kinematic viscosity $\nu$ is chosen so that
the Reynolds number
$\text{Re} = LU_{\infty}/\nu$
is either $2\times 10^4$ or $10^5$,
where $L=1$ is the side length of the unit square and
$U_{\infty} = \max(u_{\theta}) = 0.068$
is the maximal velocity magnitude.

The time span $[0, 60]$ is made long enough
for the turbulent boundary layers to
develop prominent Lagrangian coherent structures.
The tests are performed using both AMR and single-level meshes.
In the case of single-level meshes,
uniform square grids with size
$h = 2^{-9}$ and $h = 2^{-10}$ are used
for $\text{Re} = 2\times 10^4$ and $\text{Re} = 10^5$,
respectively.
In the case of AMR,
the number of refinements of the unit square $(0, 1)^2$ is allowed to
vary between 6 and 9 for $\text{Re} = 2\times 10^4$,
and between 7 and 10 for $\text{Re} = 10^5$.
Since we are only interested in regions near the domain boundaries,
the number of refinements within the disk of radius 0.2 centered at the vortex center is
limited to 7 and 8
for $\text{Re} = 2\times 10^4$ and $\text{Re} = 10^5$,
respectively.
After every 50 time steps,
we adaptively refine and coarsen the mesh and
adjust the time step size,
where the thresholds for refinement and coarsening
in the D\"orfler marking strategy (\ref{eq:dorflerStrategy})
are $\theta_{\mathrm{R}} = 0.6$ ($\text{Re} = 2\times 10^4$),
$\theta_{\mathrm{R}} = 0.8$ ($\text{Re} = 10^5$),
and $\theta_{\mathrm{C}} = 0.1$.
The time step size is adjusted by (\ref{eq:courantNumber}) and
$\text{Cr} = 0.8$.

\begin{figure}
  \centering
  \subfigure[$\text{Re} = 2\times 10^4$, single-level mesh with $h = 2^{-9}$]{
    \includegraphics[width=0.48\textwidth]{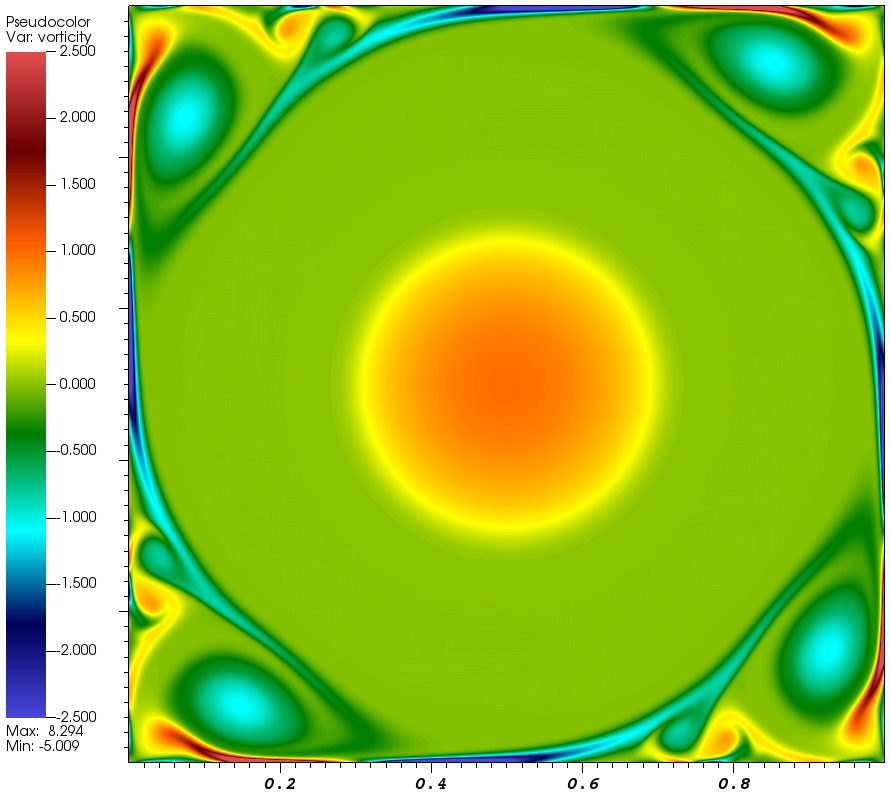}
  }
  \hfill
  \subfigure[$\text{Re} = 2\times 10^4$, an AMR hierarchy with
  \mbox{$h\in [2^{-9}, 2^{-6}]$}]{
    \includegraphics[width=0.48\textwidth]{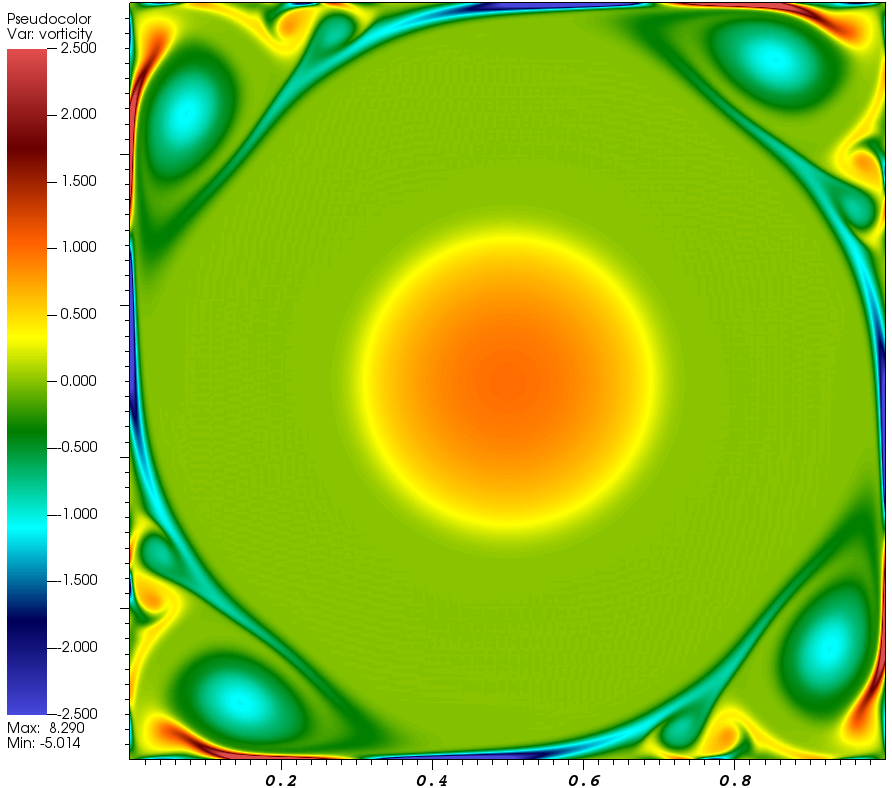}
  }

  \subfigure[$\text{Re} = 10^5$, single-level mesh with $h = 2^{-10}$]{
    \includegraphics[width=0.48\textwidth]{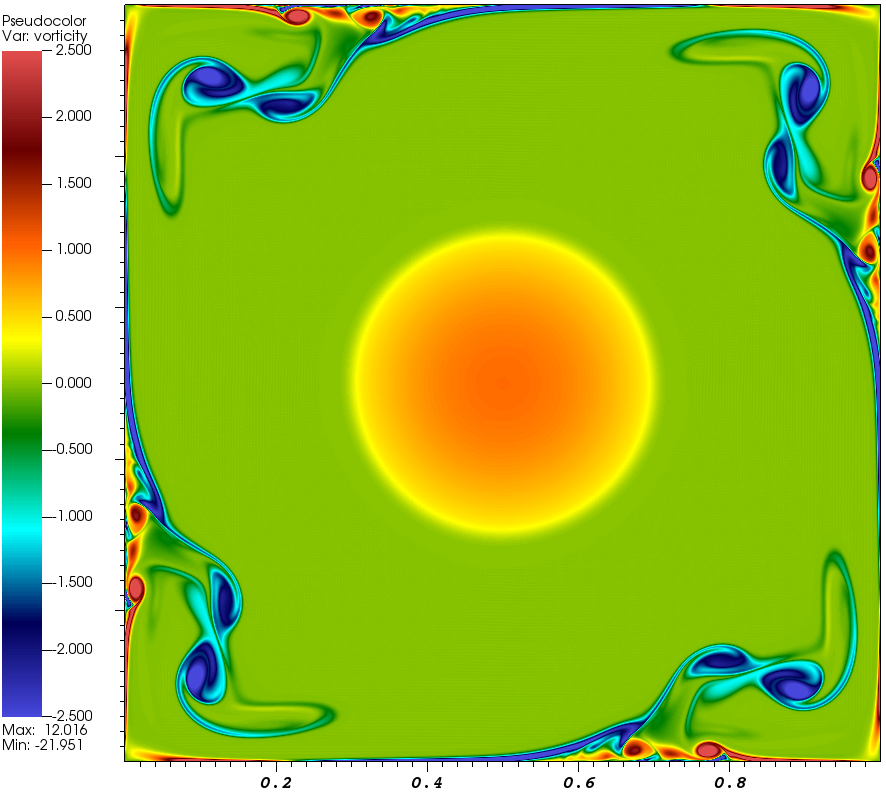}
  }
  \hfill
  \subfigure[$\text{Re} = 10^5$, an AMR hierarchy with
  \mbox{$h\in [2^{-10}, 2^{-7}]$}]{
    \includegraphics[width=0.48\textwidth]{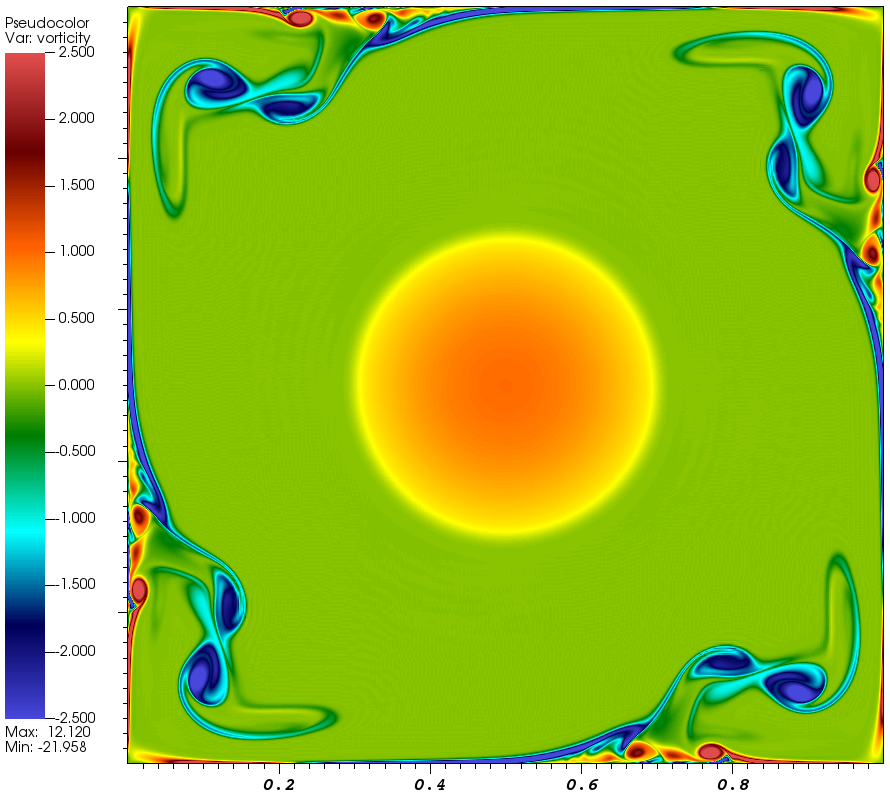}
  }
  \caption{Results of GePUP-FEM of the vorticity field
    at $t_{\mathrm{e}} = 60$
    for the single-vortex test.}
  \label{fig:singleVortexVorticity}
\end{figure}

\begin{figure}
  \centering
  \subfigure[$\text{Re} = 2\times 10^4$]{
    \includegraphics[width=0.48\textwidth]{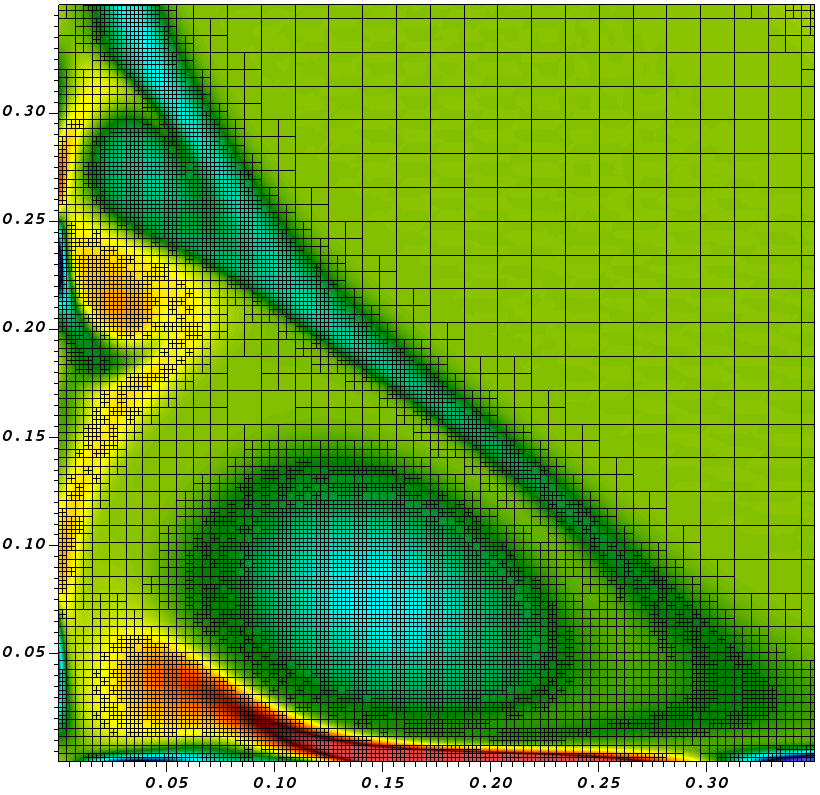}
  }
  \hfill
  \subfigure[$\text{Re} = 10^5$]{
    \includegraphics[width=0.48\textwidth]{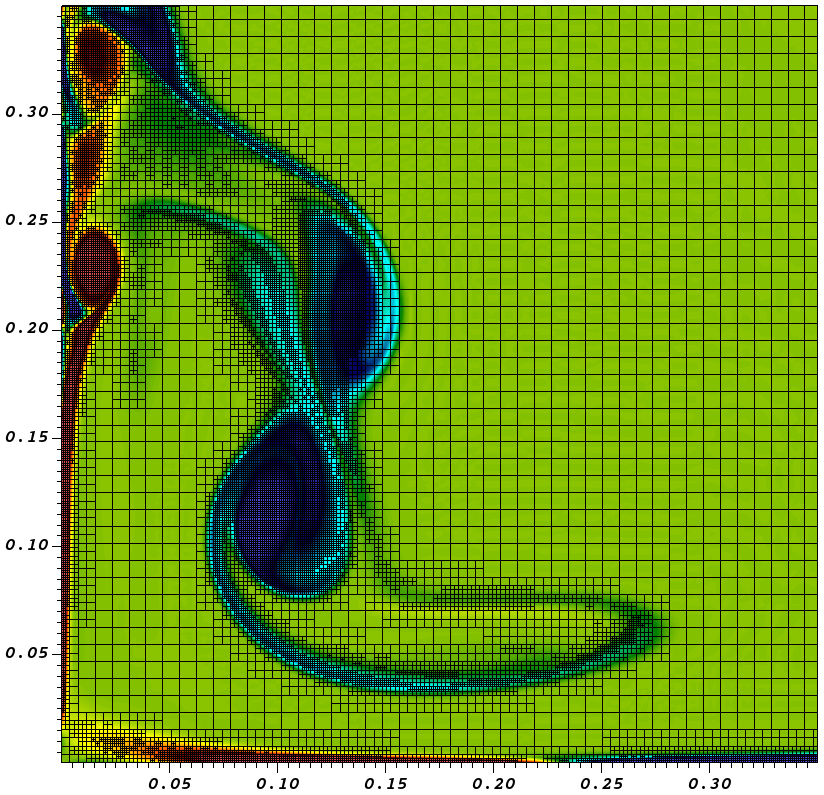}
  }
  \caption{The vorticity field of GePUP-FEM at $t_{\mathrm{e}} = 60$
    on the locally refined mesh
    within $(0, 0.35)^2$
    for the single-vortex test.}
  \label{fig:singleVortexAMRZoomIn}
\end{figure}

Figure \ref{fig:singleVortexVorticity} shows
snapshots of the vorticity field at $t_{\mathrm{e}} = 60$.
The essential features of vortex sheet roll-up and counter-vortices
agree well with those in \cite{bell91}.
In terms of the coherent structures,
there is no observable difference between the single-level snapshots
and their AMR counterparts.
A close-up view of the vorticity field and the adaptive mesh
is also shown in Figure \ref{fig:singleVortexAMRZoomIn},
demonstrating the effectiveness of our adaptive solver.
We have also implemented the UPPE formulation (\ref{eq:UPPEstrong})
and found that
the $L^2$-norm of the velocity divergence remains around $10^{-1}$,
leading to qualitatively incorrect results.

We also compare the performance of AMR to
that of a single-level mesh
in Table \ref{tab:AMRPerformanceSingleVortex},
where the four-level AMR hierarchy uses far fewer computational resources
than the single-level mesh;
for example,
at $\text{Re} = 2\times 10^4$ and $\text{Re} = 10^5$,
\revise{
  the average number of elements in AMR is only
  10.7\% and 6.36\% of that in the single-level mesh,
  respectively;
  the average number of pressure DoFs in AMR is only
  11.4\% and 6.75\% of that in the single-level mesh,
  respectively,
}
and the CPU time consumed by AMR is only 20.4\% and 9.37\% of
that consumed by the single-level mesh,
respectively.
This clearly illustrates that AMR is very promising
in simulating high-Reynolds-number flows.
However,
the relative savings of AMR on computational resources
depend drastically on problem-specific aspects
such as the refinement criteria
and the number of levels.
Hence the reader should not regard the results
in Table \ref{tab:AMRPerformanceSingleVortex}
as universal estimates
but rather just as one particular example of the efficiency of AMR.

\begin{table}
  \caption{Performance comparison
    between AMR and single-level meshes
    for the single-vortex test at $\text{Re} = 2\times 10^4$
    and $\text{Re} = 10^5$.
    The single level is discretized by square grids with size 
    $h = 2^{-9}$ ($\text{Re} = 2\times 10^4$) and
    $h = 2^{-10}$ ($\text{Re} = 10^5$).
    The AMR hierarchy contains four levels with
    $h$ varying between $2^{-9}$ and $2^{-6}$
    for $\text{Re} = 2\times 10^4$,
    and between $2^{-10}$ and $2^{-7}$
    for $\text{Re} = 10^5$.
    The simulations were performed on
    a server with AMD Ryzen Threadripper PRO 3995WX 64-Core CPUs
    running at 2.70 GHz.
    \revise{
      The numbers of elements, velocity and pressure DoFs are averaged values
      when AMR is used.
      The CPU time is the total accumulated time over the entire simulation.
    }}
  \centering
  \renewcommand{\arraystretch}{1.5}
  \begin{tabular}{c|c|c|c|c|c}
 \hline
 \multicolumn{2}{r|}{} & \revise{Elements} & Velocity DoFs & Pressure DoFs & CPU time (seconds)
 \\ \hline
 \multirow{2}*{$\text{Re} = 2\times 10^4$} & Single-level 
 & \revise{262,144} & 4,724,738 & 2,362,369 & 10,032
 \\ \cline{2-6}
 & AMR & \revise{27,940} & \revise{540,104} & \revise{270,052} & 2,049
 \\ \hline
 \multirow{2}*{$\text{Re} = 10^5$} & Single-level
 & \revise{1,048,576} & 18,886,658 & 9,443,329 & 88,604
 \\ \cline{2-6}
 & AMR & \revise{66,665} & \revise{1,275,076} & \revise{637,538} & 8,302
 \\ \hline
\end{tabular}

  \label{tab:AMRPerformanceSingleVortex}
\end{table}

\revise{
\subsubsection{Lid-driven cavity at $\mathrm{Re} = 10^4$}
\label{sec:lid-driven-cavity}

This test examines the well-known 2D lid-driven cavity flow
\cite{ghia82:_high_re_navier_stokes,erturk05:_numer_reynol},
which undergoes a Hopf bifurcation
from steady to unsteady periodic flow
at a critical Reynolds number
$\mathrm{Re}_{\mathrm{c}}\in [8017.6, 8018.8)$
\cite{auteri02:_numer}.
We focus on a high Reynolds number of $\mathrm{Re} = 10^4$,
where the flow exhibits unsteady periodic behavior
with multiple dominant frequencies \cite{bruneau06}.
The computational domain is the unit square $(0, 1)^2$,
with the top lid moving horizontally at a constant velocity $U=1$
to drive the flow,
while the other three walls remain stationary.
The simulation begins with an impulsive start,
i.e.,
$\bm{u}_0 = \bm{0}$,
and the forcing term is set to $\bm{f} = \bm{0}$.

The parameters in this test are as follows.
The time span $[0, 500]$ is chosen sufficiently long
to ensure the flow reaches a fully developed periodic state.
After every 50 time steps,
we adaptively refine and coarsen the mesh and
adjust the time step size,
where the number of refinements of the unit square $(0, 1)^2$ is allowed to
vary between 6 and 8,
and the thresholds for refinement and coarsening
in the D\"orfler marking strategy (\ref{eq:dorflerStrategy}) are
$\theta_{\mathrm{R}} = 0.6$ and $\theta_{\mathrm{C}} = 0.1$,
respectively.
The time step size is adjusted by (\ref{eq:courantNumber})
with a Courant number of $\mathrm{Cr} = 1.2$.

The results of GePUP-FEM at the final time $t_{\mathrm{e}} = 500$ are
presented in Figure \ref{fig:lidDrivenCavityStreamline}.
The velocity profiles along the cavity centerlines
(i.e., $u(0.5, y)$ and $v(x, 0.5)$)
show excellent agreement with the benchmark data from
Ghia et al. \cite{ghia82:_high_re_navier_stokes}
and Erturk et al. \cite{erturk05:_numer_reynol},
particularly in the boundary layer regions near the walls.
Table \ref{tab:lidDrivenCavityVelMinMax} lists
the maximum and minimum horizontal velocity values and
their corresponding locations along the cavity centerline.
Our results fall within the range of
those reported by Ghia et al. \cite{ghia82:_high_re_navier_stokes}
and Erturk et al. \cite{erturk05:_numer_reynol}.
This test further demonstrates the capability of GePUP-FEM
in accurately and efficiently resolving
the underlying physics of high-Reynolds-number flows.

\begin{figure}
  \centering
  \subfigure[\revise{streamlines}]{
    \includegraphics[width=0.48\textwidth]{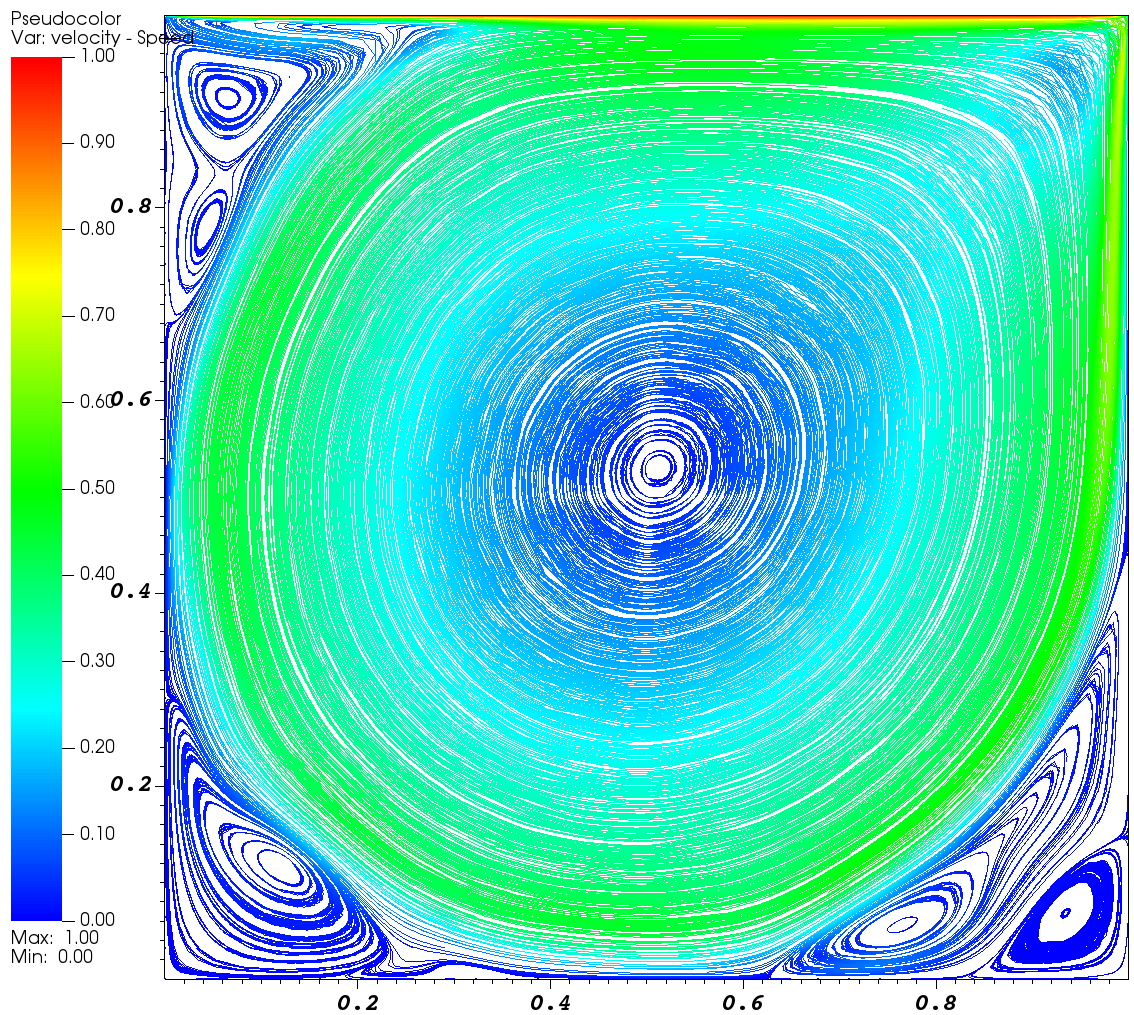}
  }
  \hfill
  \subfigure[\revise{velocity profiles along the centerlines}]{
    \includegraphics[width=0.48\textwidth]{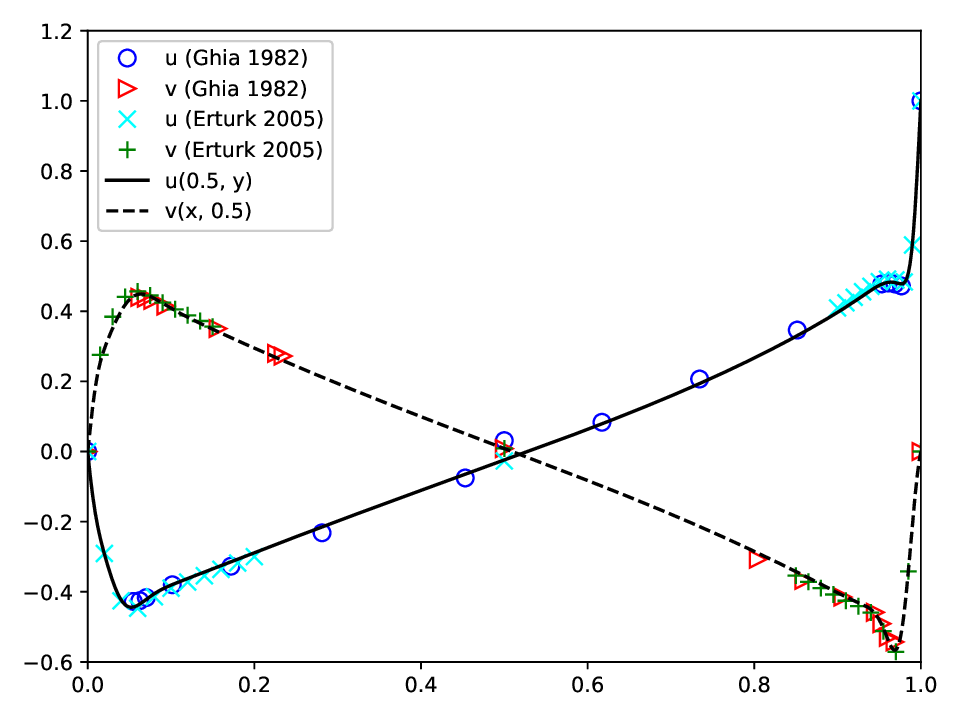}
  }
  \caption{\revise{Results of GePUP-FEM for the lid-driven cavity flow
    at $\mathrm{Re} = 10^4$ and $t = 500$.}}
  \label{fig:lidDrivenCavityStreamline}
\end{figure}

\begin{table}
  \caption{\revise{Maximum and minimum horizontal velocity values
    and their corresponding $x$-coordinates
    along the cavity centerline
    for the lid-driven cavity flow at $\mathrm{Re} = 10^4$ and $t = 500$.}}
  \centering
  \renewcommand{\arraystretch}{1.5}
  \revise{\begin{tabular}{ccccc}
  \hline
  & $x_{\max}$ & $v_{\max}$ & $x_{\min}$ & $v_{\min}$
  \\ \hline
  Ghia et al. \cite{ghia82:_high_re_navier_stokes} & 0.063 & 0.4398 & 0.970 & $-0.5430$
  \\
  Erturk et al. \cite{erturk05:_numer_reynol} & 0.060 & 0.4566 & 0.970 & $-0.5712$
  \\
  Present & 0.061 & 0.4516 & 0.969 & $-0.5664$
  \\ \hline
\end{tabular}

}
  \label{tab:lidDrivenCavityVelMinMax}
\end{table}
}

\subsubsection{Flow past a circular cylinder}
\label{sec:flow-past-cylinder}

In this test,
we consider the popular
benchmark of
laminar flows past a circular cylinder,
where the flow characteristics depend primarily on the Reynolds number.
Experimental studies and numerical simulations
\cite{williamson96:_vortex_dynam_cylin_wake}
indicate that
up to a Reynolds number
(based on the freestream velocity and the cylinder diameter)
of about 49
the flow is steady and symmetric with respect to the wake centerline.
For higher Reynolds numbers,
the flow becomes unsteady and
exhibits periodic von Karman vortex shedding,
exerting oscillatory drag and lift forces on the cylinder.
The flow remains two-dimensional up to a Reynolds number of about 190, 
beyond which the flow becomes intrinsically three-dimensional.

The problem setup consists of
a rectangular channel $(0, 32)\times (0, 16)$
containing a circular cylinder of diameter $D = 1$ centered at $(8, 8)$.
On the left, top, and bottom walls,
we impose for the velocity
the Dirichlet boundary conditions
\begin{equation*}
  \bm{u}(x, y, t) = (U_{\infty}\omega(t), 0)^T, \qquad
  U_{\infty} = 1, \quad
  \omega(t) =
  \begin{cases}
    \frac{1}{2} - \frac{1}{2}\cos(5\pi t) & t\le \frac{1}{5}, \\
    1 & t > \frac{1}{5}, 
  \end{cases}
\end{equation*}
where $\omega(t)$ regularizes the startup phase of the flow.
On the cylinder surface,
no-slip boundary conditions $\bm{u} = \bm{0}$ are enforced for the velocity,
whereas on the right wall,
homogeneous Neumann and homogeneous Dirichlet boundary conditions
are enforced for
the velocity $\bm{u}$ and the pressure $p$, respectively.
There is no external force,
i.e., $\bm{f} = \bm{0}$.
The density is set to $\rho = 1$ and
the kinematic viscosity $\nu$ is chosen so that
the Reynolds number $\text{Re} = U_{\infty}D/\nu$ is either 100 or 200.

For the right boundary of the domain,
 we impose the outflow condition
 by modifying boundary conditions
 in the GePUP formulation (\ref{eq:GePUP}) as follows: 
\begin{itemize}
\item
  for the momentum equation (\ref{eq:GePUP}a),
  the condition of $\bm{w}$ is homogeneous Neumann on the outflow boundary;
\item
  for the projection equation (\ref{eq:GePUP}c),
  we impose on $\phi$
  homogeneous Dirichlet boundary conditions on the outflow boundary
  (so that the tangential velocity is preserved)
  and homogeneous Neumann boundary conditions on other walls
  (so that the normal velocity is preserved);
\item
  for the pressure extraction equation (\ref{eq:GePUP}e),
  we impose on $q$
  homogeneous Dirichlet boundary conditions on the outflow boundary and
  the Neumann boundary condition (\ref{eq:GePUP}f) on other walls.
\end{itemize}

The parameters in this test are as follows.
The time span is $[0, 200]$.
After every 30 time steps,
we adaptively refine and coarsen the mesh and adjust the time step size,
where the thresholds for refinement and coarsening
in the D\"orfler marking strategy (\ref{eq:dorflerStrategy})
are $\theta_{\mathrm{R}} = 0.6$ (Re = 100),
$\theta_{\mathrm{R}} = 0.8$ (Re = 200),
and $\theta_{\mathrm{C}} = 0.1$.
The time step size is adjusted by (\ref{eq:courantNumber}) and Cr = 1.2.
To avoid a very small time step size caused by excessive mesh refinements,
the initial mesh is refined at most $r_{\max}$ times.

\begin{figure}
  \centering
  \subfigure[Re = 100, $t = 173.12$]{
    \includegraphics[width=0.75\textwidth]{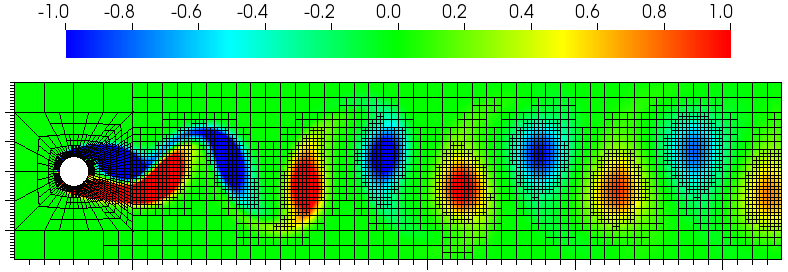}
  }

  \subfigure[Re = 200, $t = 174.13$]{
    \includegraphics[width=0.75\textwidth]{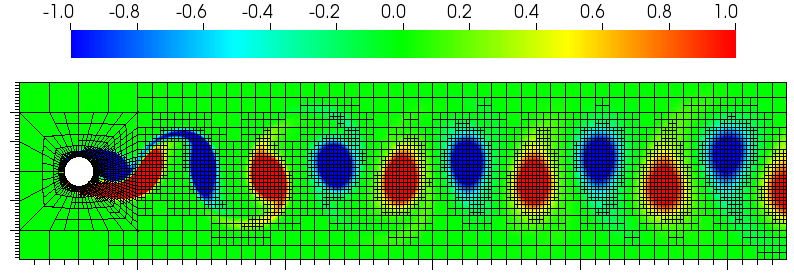}
  }
  \caption{Snapshots of
    the vorticity field on the locally refined mesh
    within the subdomain $(6, 32)\times (5, 11)$
    for the flow past a cylinder.
    The initial mesh is refined
    at most $r_{\max} = 3$ times.}
  \label{fig:flowPastACylinderVorticityAndMesh}
\end{figure}

Snapshots of the vorticity field on the locally refined mesh
within the subdomain $(6, 32)\times(5, 11)$
are shown in Figure \ref{fig:flowPastACylinderVorticityAndMesh},
which clearly reveals the von Karman vortex street phenomenon
for both Re = 100 and Re = 200.
Refining the mesh with $r_{\max} = 3$ results in
a mesh resolution of about $0.03D$ in the vicinity of the cylinder.
During the periodic shedding of vortices,
the number of pressure DoFs is
approximately $4.7\times 10^4$,
whereas refining every element in the initial mesh three times
would increase the number of pressure DoFs to 305,472!
This clearly illustrates
how AMR algorithms can substantially reduce the problem size,
thereby enhancing computational efficiency.

To validate our numerical solver
against established results in the literature,
we calculate the drag and lift coefficients,
$C_{\mathrm{D}}(t)$ and $C_{\mathrm{L}}(t)$,
obtained by evaluating
the $x$- and $y$-components of the surface integral
\begin{equation*}
  \frac{2}{\rho U_{\infty}^2D}\int_S
  (-q\bm{I} + 2\mu\bm{\sigma})\bm{n},
\end{equation*}
where $S$ is the surface of the cylinder,
$\bm{I}$ the identity tensor,
$\mu = \rho \nu$ the dynamic viscosity,
$\bm{\sigma} = \left(\nabla \bm{u} + \nabla\bm{u}^T\right)/2$
the symmetric stress tensor,
and $\bm{n}$ the unit outward normal to the cylinder.
We also calculate the dimensionless Strouhal number
$\text{St} = fD/U_{\infty}$,
where $f$ is the vortex shedding frequency 
computed from the temporal variation of
the lift coefficient $C_{\mathrm{L}}(t)$.

\begin{figure}
  \centering
  \subfigure[$C_{\mathrm{D}}(t)$ for Re = 100]{
    \includegraphics[width=0.48\textwidth]{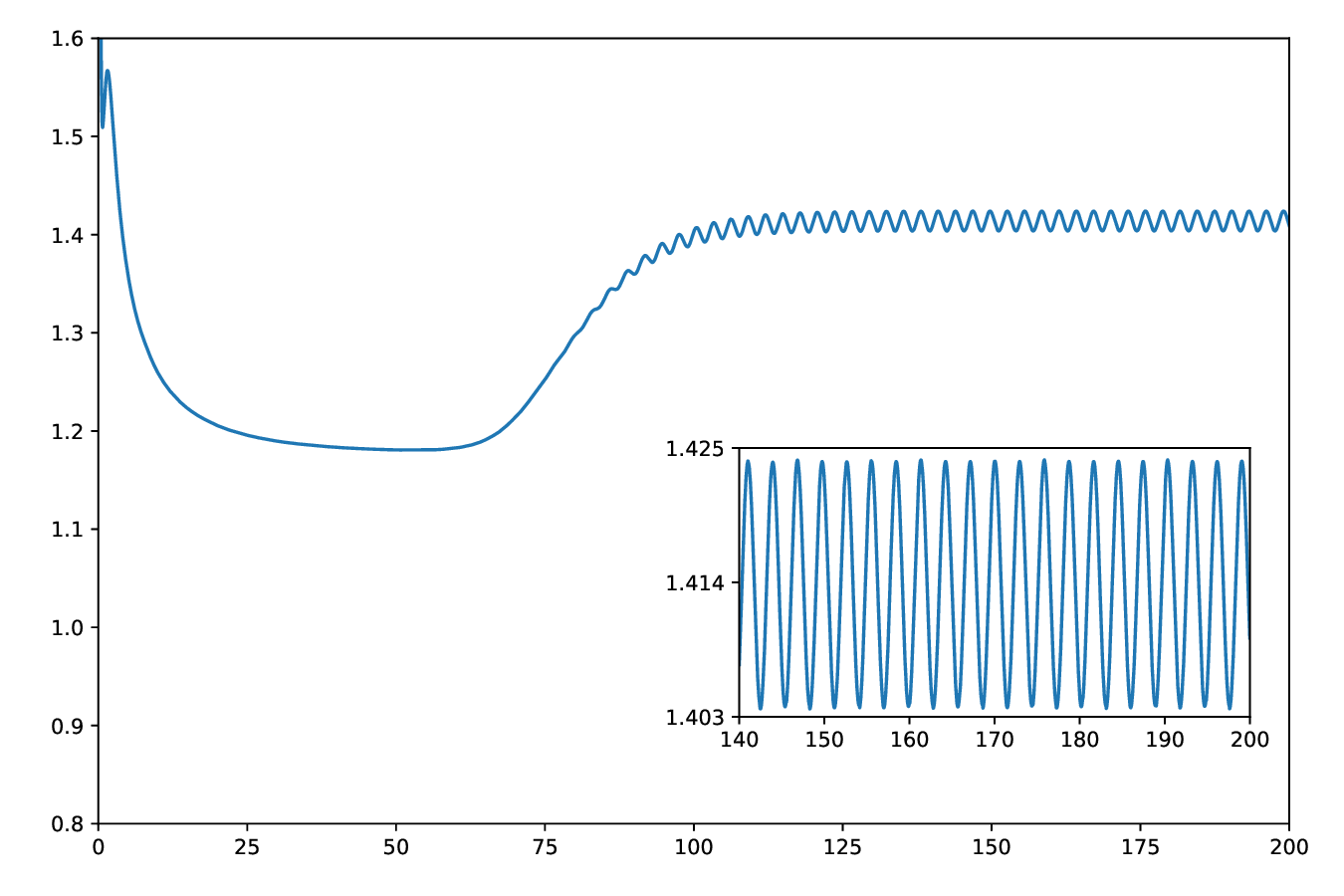}
  }
  \hfill
  \subfigure[$C_{\mathrm{D}}(t)$ for Re = 200]{
    \includegraphics[width=0.48\textwidth]{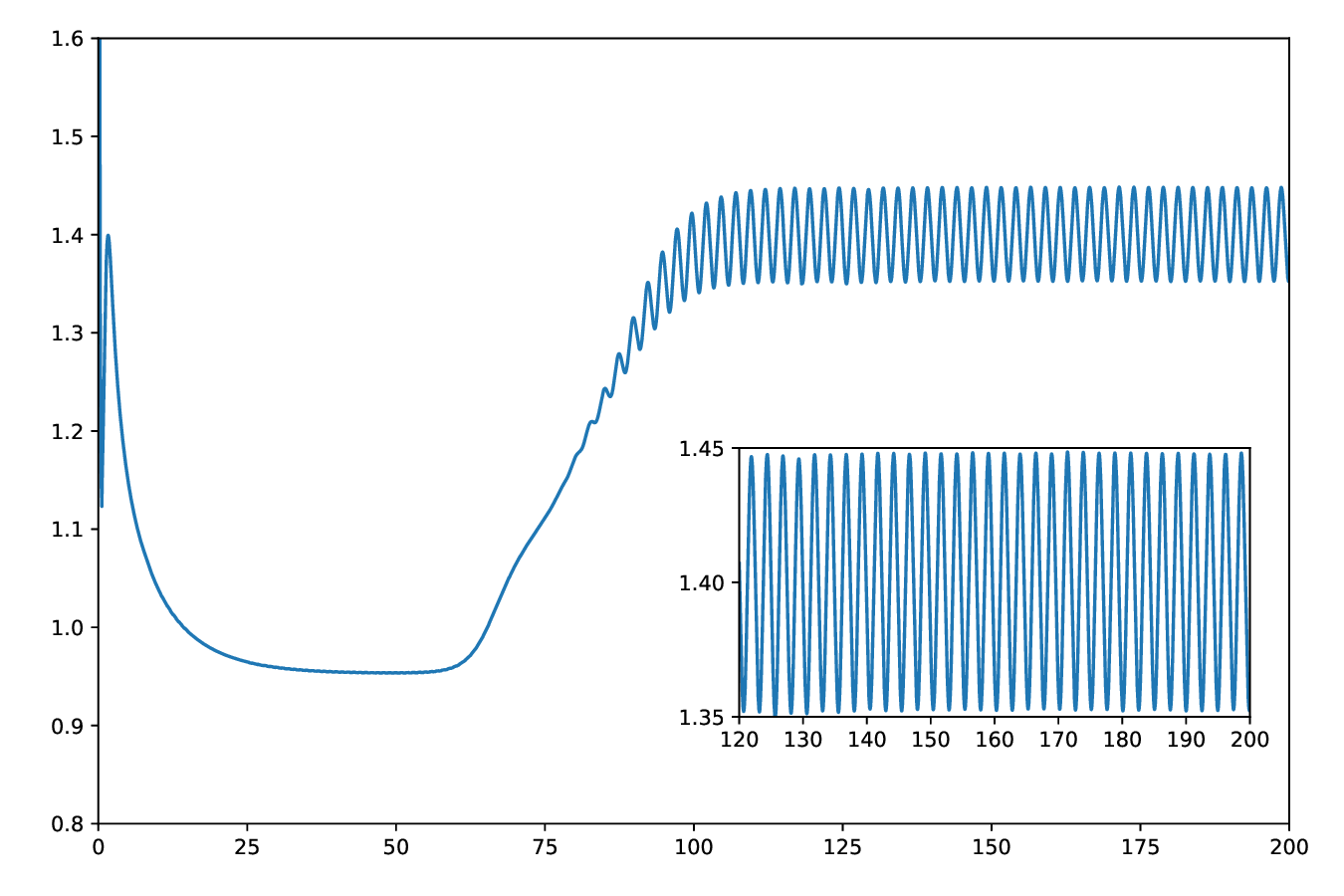}
  }
  
  \subfigure[$C_{\mathrm{L}}(t)$ for Re = 100]{
    \includegraphics[width=0.48\textwidth]{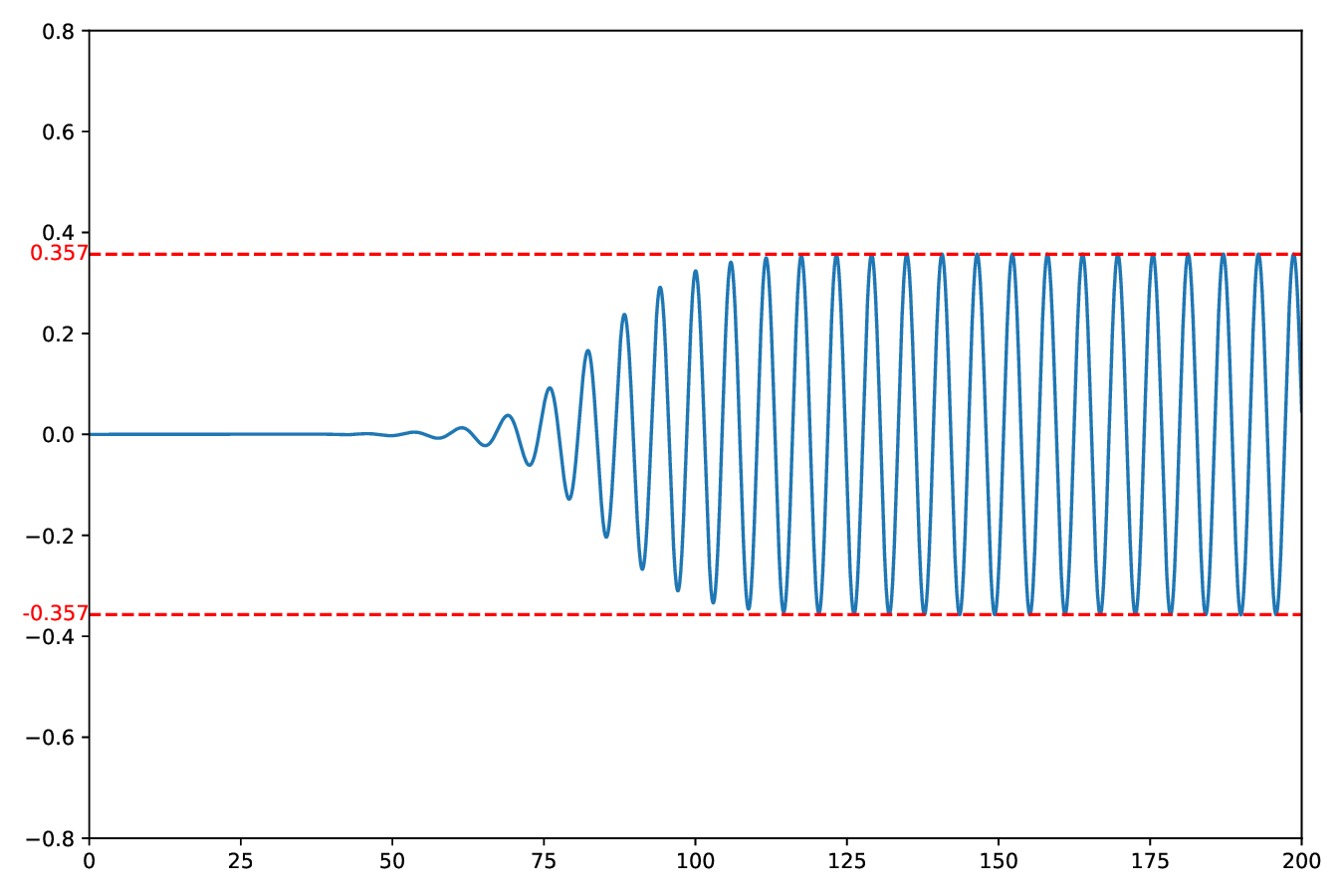}
  }
  \hfill
  \subfigure[$C_{\mathrm{L}}(t)$ for Re = 200]{
    \includegraphics[width=0.48\textwidth]{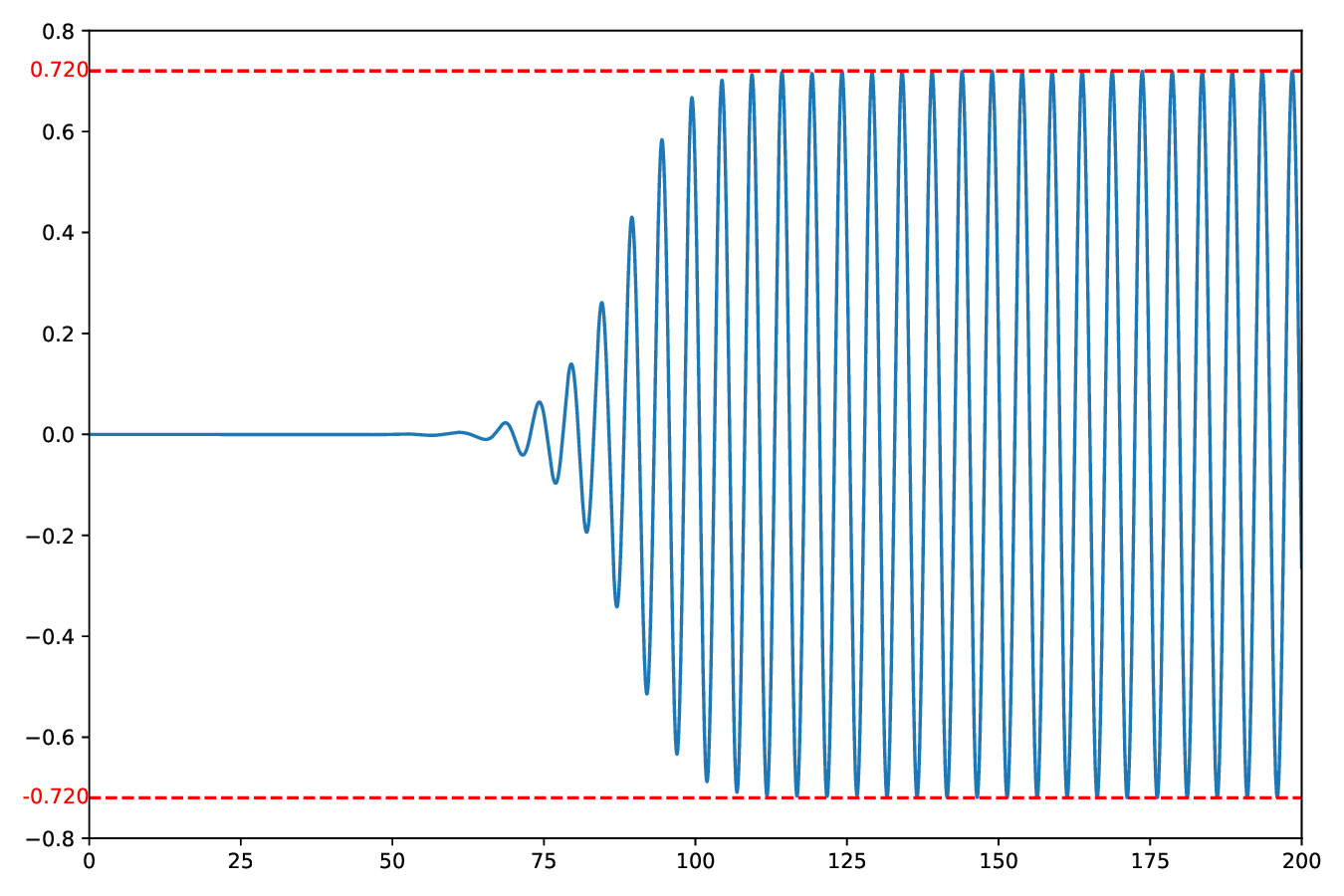}
  }
  \caption{Temporal variations of the drag and lift coefficients
    for the flow past a circular cylinder.
    The abscissa represents time
    while an ordinate is the calculated value of a coefficient.
    The initial mesh is refined
    at most $r_{\max} = 4$ times.
    The results for $r_{\max} = 3$ are similar except that
    the vortex shedding occurs at a later time.}
  \label{fig:flowPastACylinderDragAndLift}
\end{figure}

Temporal variations of the drag and lift coefficients
for Reynolds numbers
Re = 100 and Re = 200
are shown in Figure \ref{fig:flowPastACylinderDragAndLift}.
For the higher Reynolds number,
vortex shedding emerges earlier in the wake behind the cylinder and
stabilizes into the Karman vortex street sooner.
In Table \ref{tab:flowPastACylinderLiftDragStrouhal},
we list values of the drag coefficient,
lift coefficient,
and Strouhal number, 
demonstrating good agreements of our results
with those in the literature.

\begin{table}
  \caption{Drag coefficient, lift coefficient, and Strouhal number for
    the flow past a circular cylinder.
    The results for $r_{\max} = 5$ are identical to
    those for $r_{\max} = 4$ up to the displayed digits.}
  \centering
  \renewcommand{\arraystretch}{1.3}
  \begin{tabular}{ccccccc}
  \hline
  & \multicolumn{2}{c}{Drag coefficient} & \multicolumn{2}{c}{Lift coefficient} & \multicolumn{2}{c}{Strouhal number}
  \\
  \cline{2-3} \cline{4-5} \cline{6-7}
  & $\text{Re} = 100$ & $\text{Re} = 200$ & $\text{Re} = 100$ & $\text{Re} = 200$ & $\text{Re} = 100$ & $\text{Re} = 200$
  \\ \hline
  Blomquist et al. \cite{blomquist24:_stabl} & $1.387\pm 0.019$ & $1.370\pm 0.060$ & $\pm 0.346$ & $\pm 0.762$ & - & -
  \\
  Braza et al. \cite{braza86:_numer} & $1.364\pm 0.015$ & $1.400\pm 0.050$ & $\pm 0.250$ & $\pm 0.750$ & 0.16 & 0.20
  \\
  Liu et al. \cite{liu98:_precon_multig_method_unstead_incom_flows} & $1.350\pm 0.012$ & $1.310\pm 0.049$ & $\pm 0.339$ & $\pm 0.690$ & 0.165 & 0.192
  \\
  Kolahdouz et al. \cite{kolahdouz20} & $1.370\pm 0.015$ & $1.390\pm 0.060$ & $\pm 0.351$ & $\pm 0.750$ & 0.168 & 0.198
  \\
  Calhoun \cite{calhoun02:_cartes_grid_method_solvin_two} & $1.330\pm 0.014$ & $1.172\pm 0.058$ & $\pm 0.298$ & $\pm 0.668$ & 0.175 & 0.202
  \\
  Mahír et al. \cite{mahir08:_numer} & $1.368\pm 0.029$ & $1.376\pm 0.048$ & $\pm 0.343$ & $\pm 0.698$ & 0.172 & 0.192
  \\
  Xu et al. \cite{xu06} & $1.423\pm 0.013$ & $1.420\pm 0.040$ & $\pm 0.340$ & $\pm 0.660$ & 0.171 & 0.202
  \\
  Present ($r_{\max} = 3$) & $1.418\pm 0.011$ & $1.402\pm 0.048$ & $\pm 0.356$ & $\pm 0.718$ & 0.172 & 0.202
  \\
  Present ($r_{\max} = 4$) & $1.414\pm 0.011$ & $1.400\pm 0.050$ & $\pm 0.357$ & $\pm 0.720$ & 0.172 & 0.202
  \\ \hline
\end{tabular}

  \label{tab:flowPastACylinderLiftDragStrouhal}
\end{table}

\subsubsection{Flow past a sphere}
\label{sec:flow-past-sphere}

In this last test,
we consider a more challenging benchmark of
three-dimensional flows past a sphere.
Similar to the flow past a circular cylinder, 
fully three-dimensional unsteady flow fields
are generated from small perturbations
despite the symmetry of the geometric configuration.
Furthermore,
these three-dimensional flows can exhibit
even more intricate kinematic and vortical structures. 
Much like the circular cylinder case,
the flow demonstrates several distinct regimes and
the transition between those regimes depends on
the Reynolds number.
At Reynolds numbers
(based on the freestream velocity and the diameter of the sphere)
between 20 and approximately 210,
the flow remains
separated, steady, axisymmetric,
and topologically similar,
whereas the flow at Reynolds numbers between 210 and 270 is
steady but non-axisymmetric.
The onset of unsteadiness occurs
at a Reynolds number in the range of 270 to 300,
where the flow pattern exhibits periodic vortex shedding
\cite{
  johnson99:_flow_reynol,
  ormieres99:_trans_turbul_wake_spher,
  tomboulides00:_numer}.

The problem setup is similar to
that in Section \ref{sec:flow-past-cylinder}. 
We consider a static sphere of diameter $D = 1$
centered at $(8, 8, 8)$
within the rectangular box
$(0, 32)\times (0, 16)\times (0, 16)$.
On the left, top, bottom, front, and back walls,
we impose for the velocity the Dirichlet boundary conditions
\begin{equation*}
  \bm{u}(x, y, z, t) = \left( U_{\infty}\omega(t), 0, 0 \right)^T, \qquad
  U_{\infty} = 1, \quad
  \omega(t) =
  \begin{cases}
    \frac{1}{2} - \frac{1}{2}\cos(5\pi t) & t\le \frac{1}{5}, \\
    1 & t > \frac{1}{5}.
  \end{cases}
\end{equation*}
On the sphere surface,
no-slip boundary conditions $\bm{u} = \bm{0}$
are enforced for the velocity,
whereas on the right wall,
homogeneous Neumann and homogeneous Dirichlet boundary conditions
are enforced for the velocity $\bm{u}$ and the pressure $p$,
respectively.
There is no external force,
i.e.,
$\bm{f} = \bm{0}$.
The density is set to $\rho = 1$ and
the kinematic viscosity $\nu$ is deduced from the Reynolds number
$\text{Re} = U_{\infty}D/\nu$,
which,
in this test,
varies from 100 to 500,
covering the different laminar flow regimes outlined above.
The outflow boundary condition is implemented
 in the same manner as that in 
 the third paragraph of Section \ref{sec:flow-past-cylinder}.

The parameters for this test are as follows.
The time span is $[0, 400]$.
After every 50 time steps,
we adaptively refine and coarsen the mesh and adjust the time step size,
where the thresholds for refinement and coarsening
in the D\"orfler marking strategy (\ref{eq:dorflerStrategy})
are $\theta_{\mathrm{R}} = 0.6$ (Re = 300),
$\theta_{\mathrm{R}} = 0.7$ (Re = 350),
$\theta_{\mathrm{R}} = 0.8$ (Re = 500),
and $\theta_{\mathrm{C}} = 0.1$.
The time step size is adjusted by
(\ref{eq:courantNumber}) and Cr = 0.8.
The initial mesh is refined at most $r_{\max} = 3$ times.
The resulting mesh resolution in the vicinity of the sphere is
approximately $0.078D$.
For Re = 500,
the number of DoFs for the pressure and the velocity is
approximately $3\times 10^6$ and $9\times 10^6$,
respectively,
which correspond to 2.62\% of the number of DoFs
for a uniform refinement of the initial mesh by three times.

To visualize vortical structures, 
we employ the Q-criterion
\cite{hunt88:_eddies_stream_conver_zones_turbul_flows}
defined by
\begin{equation}
  \label{eq:QCriterion}
  Q := \frac{1}{2}\left( \|\bm{\Omega}\|_2^2 - \|\bm{\sigma}\|_2^2 \right),
\end{equation}
where $\bm{\Omega} = \left( \nabla\bm{u} - \nabla\bm{u}^T \right) / 2$
is the vorticity tensor,
$\bm{\sigma}$ is the symmetric stress tensor,
and $\|\cdot\|_2$ denotes the Frobenius norm.
A vortex can be identified as the region surrounded by
an isosurface of $Q > 0$,
i.e.,
a region where the Frobenius norm of the vorticity tensor dominates
that of the symmetric stress tensor.

\begin{figure}
  \centering
  \subfigure[Re = 300, $t = 349.27$]{
    \includegraphics[width=0.66\textwidth]{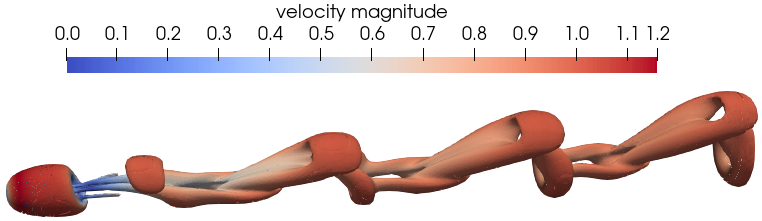}
  }

  \subfigure[Re = 350, $t = 375.15$]{
    \includegraphics[width=0.66\textwidth]{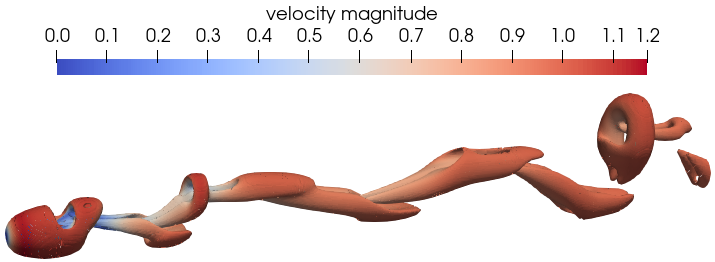}
  }
  \caption{
    Instantaneous vortical structures
    identified by the Q-criterion (\ref{eq:QCriterion})
    for the unsteady flow past a sphere.
    Isocontours of the Q-criterion ($Q = 0.001$)
    colored by the velocity magnitude
    are shown.}
  \label{fig:flowPastASphereQCriterion}
\end{figure}

Snapshots of the vortical structures for
the unsteady flows at Re = 300 and Re = 350
are shown in Figure \ref{fig:flowPastASphereQCriterion}, 
where one clearly observes periodical shedding of
harpin vortices from the sphere,
the ubiquity of vortex loops in the wake,
and the spanwise symmetry with respect to a vertical plane.
These results are consistent with those reported in
\cite{johnson99:_flow_reynol, tomboulides00:_numer},
implying that the three-dimensional vortical structures
have been successfully captured by GePUP-FEM.

As in the previous test,
we also evaluate the nondimensional force over the static sphere
\begin{equation*}
  \bm{F} = \frac{8}{\rho U_{\infty}^2\pi D^2}\int_S
  ( -q\bm{I} + 2\mu \bm{\sigma} )\bm{n},
\end{equation*}
where $S$ is the surface of the sphere and
$\bm{n}$ the unit outward normal to the sphere.
The drag coefficient $C_{\mathrm{D}}(t)$ and
the lift coefficients
$C_{\mathrm{L}}(t)$, $C_{\mathrm{L}, y}(t)$, $C_{\mathrm{L}, z}(t)$
are defined by
\begin{equation*}
  C_{\mathrm{D}}(t) = \bm{F}\cdot\bm{e}_x,
  \quad
  C_{\mathrm{L}}(t) = \sqrt{C_{\mathrm{L}, y}(t)^2 + C_{\mathrm{L}, z}(t)^2},
  \quad
  C_{\mathrm{L}, y}(t) = \bm{F}\cdot\bm{e}_y,
  \quad
  C_{\mathrm{L}, z}(t) = \bm{F}\cdot\bm{e}_z.
\end{equation*}

\begin{figure}
  \centering
  \subfigure[$C_{\mathrm{D}}(t)$]{
    \includegraphics[width=0.48\textwidth]{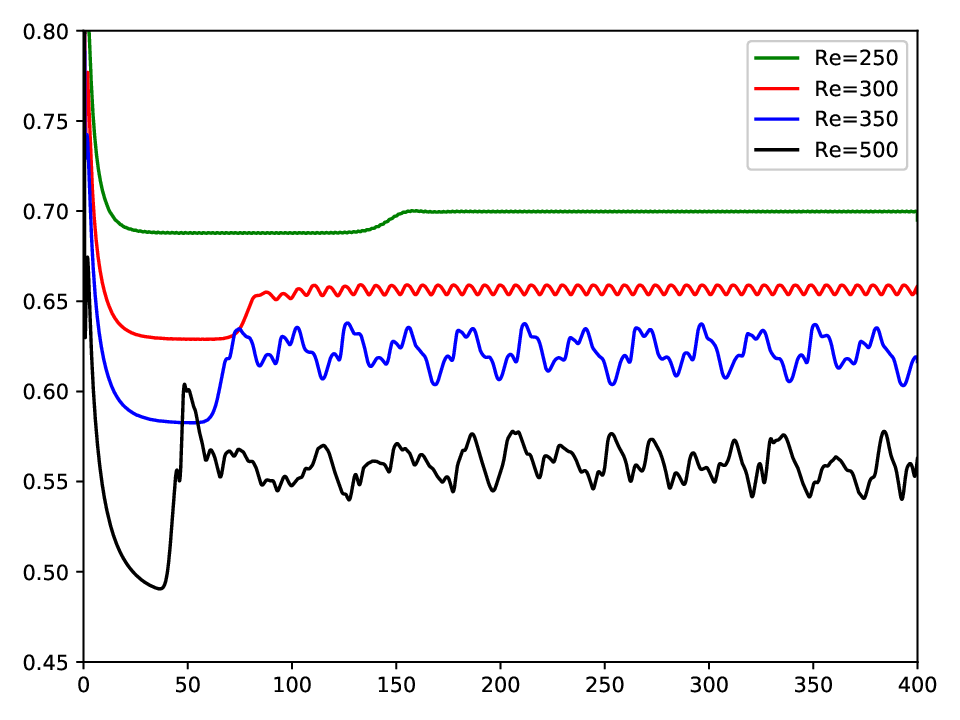}
  }

  \subfigure[$C_{\mathrm{L}, y}(t)$]{
    \includegraphics[width=0.48\textwidth]{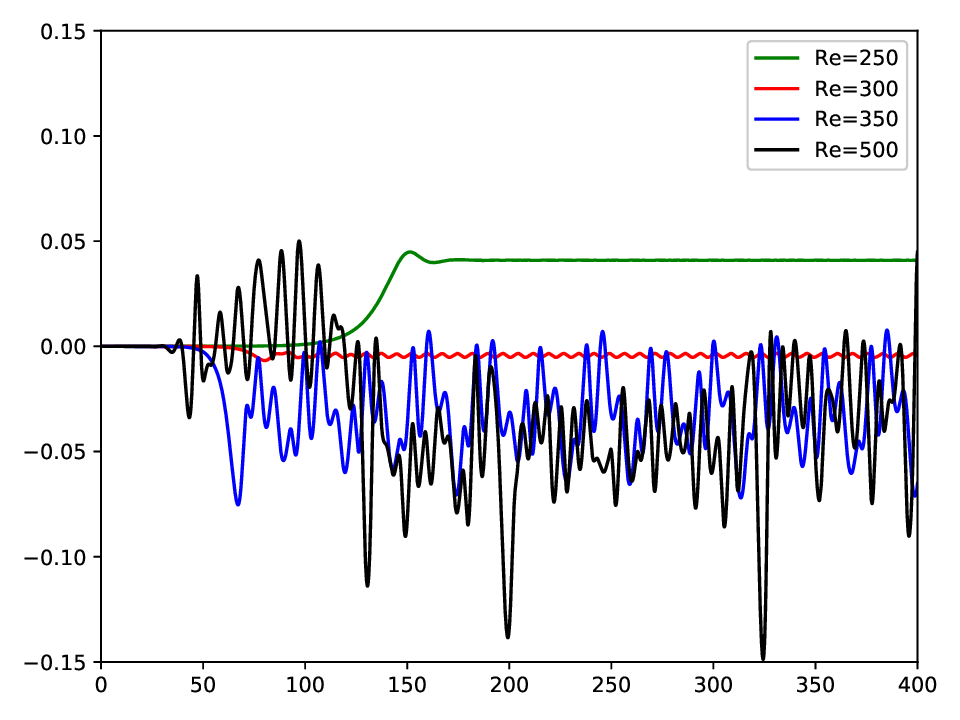}
  }
  \hfill
  \subfigure[$C_{\mathrm{L}, z}(t)$]{
    \includegraphics[width=0.48\textwidth]{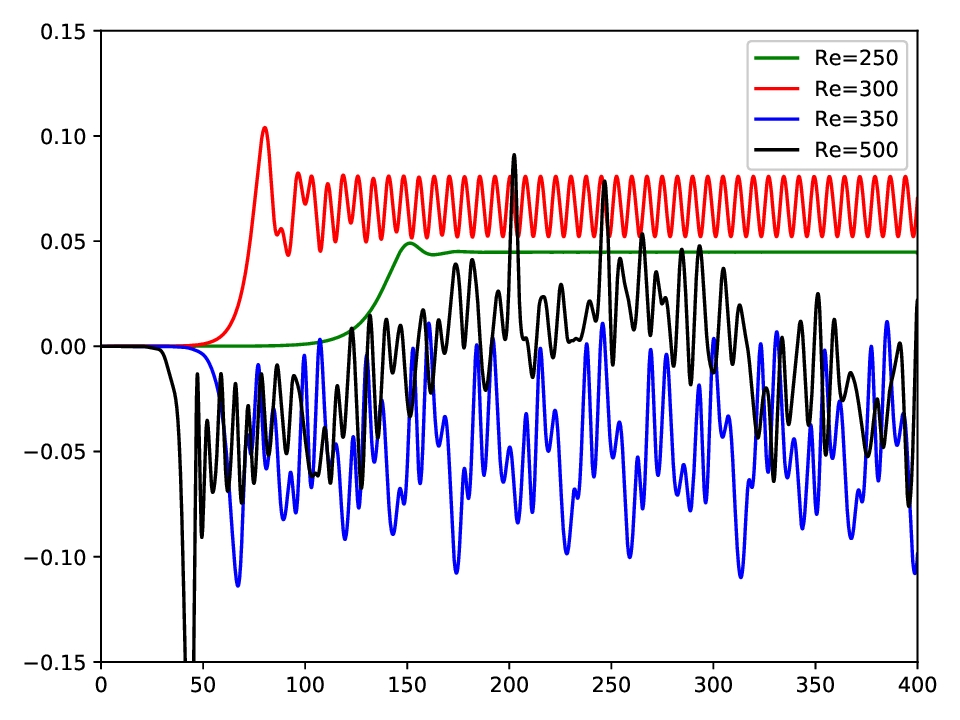}
  }
  \caption{Temporal variations of the drag and lift coefficients
    for the flow past a sphere.
    The abscissa represents time
    while an ordinate is the calculated value of a coefficient.
  }
  \label{fig:flowPastASphereDragLiftCoefs}
\end{figure}

Due to the axisymmetry of the flow at Re = 100 and Re = 200,
the lift coefficient $C_{\mathrm{L}}(t)$ should be identically zero
\cite{johnson99:_flow_reynol}.
Our numerical results are $3.86\times 10^{-13}$ and $8.91\times 10^{-5}$,
both being close to zero.
The temporal variations of the drag coefficient $C_{\mathrm{D}}(t)$
and the lift coefficients $C_{\mathrm{L}, y}(t)$ and $C_{\mathrm{L}, z}(t)$
at Reynolds numbers 250, 300, 350, and 500 are shown in
Figure \ref{fig:flowPastASphereDragLiftCoefs}.
The lift coefficients $C_{\mathrm{L}, y}(t)$ and $C_{\mathrm{L}, z}(t)$
at Re = 250
approach nonzero constants
after a nondimensional time of approximately $tU_{\infty}/D = 150$;
this is consistent with the experimental observation that
the flow is non-axisymmetric
but remains steady.
At Re = 300 and Re = 350,
the flow is unsteady and
vortices periodically shed from the static sphere
\cite{bagchi00:_direc_numer_simul_flow_heat, mittal08}.
We evaluate the Strouhal number
$\text{St} = fD/U_{\infty}$,
where the vortex shedding frequency $f$ is obtained by
calculating the averaged period between successive peak values
in the temporal variation of the lift coefficient $C_{\mathrm{L}, z}(t)$.

Our numerical results of the Strouhal number
and 
the averaged drag coefficient $C_{\mathrm{D}}(t)$
and the averaged lift coefficient $C_{\mathrm{L}}(t)$
within the time interval $[200,400]$
are summarized and compared with available data from the literature
in Tables \ref{tab:flowPastASphereAverageDrag} and
\ref{tab:flowPastASphereAverageLiftStrouhal}.
Our results agree well with those in previous studies,
confirming the accuracy of GePUP-FEM
over a relatively wide range of Reynolds numbers.

\begin{table}
  \caption{Values of the averaged drag coefficient
    within $t\in [200, 400]$
    for the flow past a sphere.}
  \centering
  \renewcommand{\arraystretch}{1.3}
  \begin{tabular}{cccccccc}
  \hline
  & $\text{Re} = 100$ & $\text{Re} = 200$ & $\text{Re} = 250$ & $\text{Re} = 300$ & $\text{Re} = 350$ & $\text{Re} = 500$
  \\ \hline
  Blomquist et al. \cite{blomquist24:_stabl} & 1.058 & 0.768 & - & 0.646 & 0.601 & -
  \\
  Johnson et al. \cite{johnson99:_flow_reynol} & 1.08 & 0.78 & 0.70 & 0.656 & - & -
  \\
  Bagchi et al. \cite{bagchi00:_direc_numer_simul_flow_heat} & 1.09 & - & 0.70 & - & 0.62 & 0.555
  \\
  Mittal et al. \cite{mittal08} & 1.08 & - & - & 0.67 & 0.625 & -
  \\
  Constantinescu et al. \cite{constantinescu03:_les_des_inves_turbul_flow} & - & - & 0.70 & 0.655 & - & -
  \\
  Kim et al. \cite{kim01:_immer_bound_finit_volum_method} & 1.087 & - & 0.701 & 0.657 & - & -
  \\
  Hartmann et al. \cite{hartmann11:_cartes} & 1.083 & 0.764 & 0.698 & 0.657 & - & -
  \\
  Giacobello et al. \cite{giacobello09:_wake_reynol} & - & - & 0.702 & 0.658 & - & -
  \\
  Present ($r_{\max} = 2$) & 1.071 & 0.765 & 0.699 & 0.657 & 0.622 & 0.558
  \\
  Present ($r_{\max} = 3$) & 1.080 & 0.767 & 0.700 & 0.657 & 0.622 & 0.560
  \\ \hline
\end{tabular}

  \label{tab:flowPastASphereAverageDrag}
\end{table}

\begin{table}
  \caption{Values of the averaged lift coefficient
    within $t\in [200, 400]$
    and the Strouhal number for
    the flow past a sphere.}
  \centering
  \renewcommand{\arraystretch}{1.3}
  \begin{tabular}{ccccc}
  \hline
  & \multicolumn{2}{c}{Lift coefficient} & \multicolumn{2}{c}{Strouhal number}
  \\
  \cline{2-3} \cline{4-5}
  & $\text{Re} = 250$ & $\text{Re} = 300$ & $\text{Re} = 300$ & $\text{Re} = 350$
  \\ \hline
  Johnson et al. \cite{johnson99:_flow_reynol} & 0.062 & 0.069 & 0.137 & -
  \\
  Bagchi et al. \cite{bagchi00:_direc_numer_simul_flow_heat} & - & - & - & 0.135
  \\
  Mittal et al. \cite{mittal08} & - & - & 0.135 & 0.142
  \\
  Constantinescu et al. \cite{constantinescu03:_les_des_inves_turbul_flow} & 0.062 & 0.065 & 0.136 & -
  \\
  Kim et al. \cite{kim01:_immer_bound_finit_volum_method} & 0.059 & 0.067 & 0.134 & -
  \\
  Hartmann et al. \cite{hartmann11:_cartes} & 0.065 & 0.069 & 0.135 & -
  \\
  Giacobello et al. \cite{giacobello09:_wake_reynol} & 0.061 & 0.067 & 0.134 & -
  \\
  Present ($r_{\max}$ = 2) & 0.060 & 0.066 & 0.134 & 0.136
  \\
  Present ($r_{\max}$ = 3) & 0.061 & 0.067 & 0.134 & 0.133
  \\ \hline
\end{tabular}

  \label{tab:flowPastASphereAverageLiftStrouhal}
\end{table}

\section{Conclusion}
\label{sec:conclusion}

We have presented GePUP-FEM,
high-order finite element methods
for solving the incompressible Navier-Stokes equations.
The weak form of the GePUP formulation is first derived and
high-order GePUP-FEM schemes with method-of-lines
are then proposed.
Results of standard benchmark problems show that
GePUP-FEM not only achieves high-order accuracy both in time and in space,
but is also capable of
accurately and efficiently resolving the right physics.
Furthermore,
GePUP-FEM offers the flexibility of choosing
finite element spaces for the velocity and the pressure 
so that the algorithmic steps
can be free of the standard inf-sup condition.

A number of future research prospects follow.
\revise{
To further enhance the stability and robustness of GePUP-FEM
in simulating flows with high Reynolds numbers,
we can employ some stabilization techniques,
such as continuous interior penalty
\cite{burman23:_implic_oseen_reynol},
where a penalty term on
the jump of the gradient over element boundaries is added.
}
To better guide the adaptive mesh refinement process,
rigorous a posteriori error estimates
need to be established
\cite{bansch16:_poster_error_estim_press_correc_schem}.
We also plan to generalize GePUP-FEM to
simulate flows with moving boundaries,
such as those in
\cite{gross11:_numer_method_two_incom_flows}
and \cite{richter17:_fluid_inter},
via incorporating high-order interface tracking methods
\cite{zhang16:_mars, zhang18:_cubic_mars_method_fourt_sixth}
and unfitted Eulerian finite element methods
\cite{burman15:_cutfem, lehrenfeld19:_eulerian_fem}.

\section*{Acknowledgments}

\revise{
We are grateful to two anonymous referees
for their insightful comments and valuable suggestions.
}
This work was supported by the grant 12272346 from
the National Natural Science Foundation of China.

\bibliographystyle{elsarticle-num}

\bibliography{GePUP-FEM.bib}

\begin{thebibliography}{10}
\expandafter\ifx\csname url\endcsname\relax
  \def\url#1{\texttt{#1}}\fi
\expandafter\ifx\csname urlprefix\endcsname\relax\def\urlprefix{URL }\fi
\expandafter\ifx\csname href\endcsname\relax
  \def\href#1#2{#2} \def\path#1{#1}\fi

\bibitem{fefferman06:_exist_navier_stokes}
C.~L. Fefferman, Existence and smoothness of the {Navier-Stokes} equation, in: J.~Carlson, A.~Jaffe, A.~Wiles (Eds.), The Millennium Prize Problems, American Mathematical Society, Providence, Rhode Island, 2006, pp. 57--67.

\bibitem{smale98:_mathem}
S.~Smale, Mathematical problems for the next century, Math. Intelligencer 20~(2) (1998) 7--15.

\bibitem{devlin03}
K.~Devlin, The Millennium Problems: The Seven Greatest Unsolved Mathematical Puzzles of Our Time, Basic Books, New York, 2003.

\bibitem{zhang16:_GePUP}
Q.~Zhang, {GePUP}: Generic projection and unconstrained {PPE} for fourth-order solutions of the incompressible {Navier-Stokes} equations with no-slip boundary conditions, J. Sci. Comput. 67~(3) (2016) 1134--1180.

\bibitem{boffi13:_mixed_finit_elemen_method_applic}
D.~Boffi, F.~Brezzi, M.~Fortin, Mixed Finite Element Methods and Applications, Springer-Verlag, Berlin Heidelberg, 2013.

\bibitem{burman07:_contin_navier_stokes}
E.~Burman, M.~A. Fern\'andez, Continuous interior penalty finite element method for the time-dependent {Navier-Stokes} equations: space discretization and convergence, Numer. Math. 107 (2007) 39--77.

\bibitem{codina07:_time}
R.~Codina, J.~Principe, O.~Guasch, S.~Badia, Time dependent subscales in the stabilized finite element approximation of incompressible flow problems, Comput. Methods Appl. Mech. Engrg. 196~(21) (2007) 2413--2430.

\bibitem{benzi05:_numer}
M.~Benzi, G.~H. Golub, J.~Liesen, Numerical solution of saddle point problems, Acta Numer. 14 (2005) 1--137.

\bibitem{chorin68:_numer_solut_navier_stokes_equat}
A.~J. Chorin, Numerical solution of the {N}avier-{S}tokes equations, Math. Comp. 22~(104) (1968) 745--762.

\bibitem{temam69:_sur_navier_stokes}
R.~Temam, Sur l’approximation de la solution des équations de {Navier–Stokes} par la méthode des pas fractionnaires {II}, Arch. Ration. Mech. Anal. 33 (1969) 377--385.

\bibitem{brown01:_accurate_projection_methods_for_ins}
D.~L. Brown, R.~Cortez, M.~L. Minion, Accurate projection methods for the incompressible {N}avier-{S}tokes equations, J. Comput. Phys. 168~(2) (2001) 464--499.

\bibitem{guermond06:_overview_of_projection_methods_for_incompressible_flows}
J.~L. Guermond, P.~Minev, J.~Shen, An overview of projection methods for incompressible flows, Comput. Methods Appl. Mech. Engrg. 195~(44-47) (2006) 6011--6045.

\bibitem{howell97}
L.~H. Howell, J.~B. Bell, An adaptive mesh projection method for viscous incompressible flow, SIAM J. Sci. Comput. 18~(4) (1997) 996--1013.

\bibitem{martin08:_navier_stokes}
D.~F. Martin, P.~Colella, D.~Graves, A cell-centered adaptive projection method for the incompressible {Navier-Stokes} equations in three dimensions, J. Comput. Phys. 227~(3) (2008) 1863--1886.

\bibitem{trebotich15:_navier_stokes}
D.~Trebotich, D.~T. Graves, An adaptive finite volume method for the incompressible {Navier-Stokes} equations in complex geometries, Commun. Appl. Math. Comput. Sci. 10~(1) (2015) 43--82.

\bibitem{blomquist24:_stabl}
M.~Blomquist, S.~R. West, A.~L. Binswanger, M.~Theillard, Stable nodal projection method on octree grids, J. Comput. Phys. 499 (2024) 112695.

\bibitem{bell89:_navier_stokes}
J.~B. Bell, P.~Colella, H.~M. Glaz, A second-order projection method for the incompressible {Navier-Stokes} equations, J. Comput. Phys. 85~(2) (1989) 257--283.

\bibitem{kim85:_applic_navier_stokes}
J.~Kim, P.~Moin, Application of a fractional-step method to incompressible {Navier-Stokes} equations, J. Comput. Phys. 59~(2) (1985) 308--323.

\bibitem{gresho87:_navier_stokes}
P.~M. Gresho, R.~L. Sani, On presssure boundary conditions for the incompressible {Navier-Stokes} equations, Int. J. Numer. Methods Fluids 7~(10) (1987) 1111--1145.

\bibitem{sanderse12:_accur_runge_kutta_navier_stokes}
B.~Sanderse, B.~Koren, Accuracy analysis of explicit {R}unge-{K}utta methods applied to the incompressible {N}avier-{S}tokes equations, J. Comput. Phys. 231~(8) (2012) 3041--3063.

\bibitem{johnston04:_accur_navier_stokes}
H.~Johnston, J.-G. Liu, Accurate, stable and efficient {N}avier-{S}tokes solvers based on explicit treatment of the pressure term, J. Comput. Phys. 199~(1) (2004) 221--259.

\bibitem{liu07:_stabil_conver_effic_navier_stokes}
J.-G. Liu, J.~Liu, R.~L. Pego, Stability and convergence of efficient {Navier-Stokes} solvers via a commutator estimate, Comm. Pure Appl. Math. 60~(10) (2007) 1443--1487.

\bibitem{liu09:_error_navier_stokes}
J.-G. Liu, J.~Liu, R.~L. Pego, Error estimates for finite-element {Navier-Stokes} solvers without standard inf-sup conditions, Chinese Ann. Math. Ser. B 30~(6) (2009) 743--768.

\bibitem{liu09:_open_navier_stokes}
J.~Liu, Open and traction boundary conditions for the incompressible {Navier-Stokes} equations, J. Comput. Phys. 228~(19) (2009) 7250--7267.

\bibitem{liu10:_stable_accurate_pressure_unsteady_incompressible_viscous_flow}
J.-G. Liu, J.~Liu, R.~L. Pego, Stable and accurate pressure approximation for unsteady incompressible viscous flow, J. Comput. Phys. 229~(9) (2010) 3428--3453.

\bibitem{jia11:_stabl_navier_stokes}
J.~Jia, J.~Liu, Stable and spectrally accurate schemes for the {Navier-Stokes} equations, SIAM J. Sci. Comput. 33~(5) (2011) 2421--2439.

\bibitem{shirokoff11:_navier_stokes}
D.~Shirokoff, R.~R. Rosales, An efficient method for the incompressible {Navier-Stokes} equations on irregular domains with no-slip boundary conditions, high order up to the boundary, J. Comput. Phys. 230~(23) (2011) 8619--8646.

\bibitem{rosales21:_high_poiss_navier_stokes}
R.~R. Rosales, B.~Seibold, D.~Shirokoff, D.~Zhou, High-order finite element methods for a pressure {Poisson} equation reformulation of the {Navier-Stokes} equations with electric boundary conditions, Comput. Methods Appl. Mech. Engrg. 373 (2021) 113451.

\bibitem{karniadakis91:_high_navier_stokes}
G.~E. Karniadakis, M.~Israeli, S.~A. Orszag, High-order splitting methods for the incompressible {Navier-Stokes} equations, J. Comput. Phys. 97~(2) (1991) 414--443.

\bibitem{petersson01:_stabil_stokes_navier_stokes}
N.~A. Petersson, Stability of pressure boundary conditions for {Stokes} and {Navier-Stokes} equations, J. Comput. Phys. 172~(1) (2001) 40--70.

\bibitem{leriche06:_numer_stokes}
E.~Leriche, E.~Perchat, G.~Labrosse, M.~O. Deville, Numerical evaluation of the accuracy and stability properties of high-order direct {Stokes} solvers with or without temporal splitting, J. Sci. Comput. 26~(1) (2006) 25--43.

\bibitem{li20:_navier_stokes}
L.~Li, A split-step finite-element method for incompressible {Navier-Stokes} equations with high-order accuracy up-to the boundary, J. Comput. Phys. 408 (2020) 109274.

\bibitem{meng20:_fourt_imex_navier_stokes}
F.~Meng, J.~W. Banks, W.~D. Henshaw, D.~W. Schwendeman, Fourth-order accurate fractional-step {IMEX} schemes for the incompressible {Navier-Stokes} equations on moving overlapping grids, Comput. Methods Appl. Mech. Engrg. 366 (2020) 113040.

\bibitem{pacheco21:_newton}
D.~R.~Q. Pacheco, R.~Schussnig, T.-P. Fries, An efficient split-step framework for {non-Newtonian} incompressible flow problems with consistent pressure boundary conditions, Comput. Methods Appl. Mech. Engrg. 382 (2021) 113888.

\bibitem{ascher97:_implic_runge_kutta}
U.~M. Ascher, S.~J. Ruuth, R.~J. Spiteri, Implicit-explicit {R}unge-{K}utta methods for time-dependent partial differential equations, Appl. Numer. Math. 25~(2) (1997) 151--167.

\bibitem{kennedy03:_addit_runge_kutta_schem_for}
C.~A. Kennedy, M.~H. Carpenter, Additive {R}unge-{K}utta schemes for convection-diffusion-reaction equations, Appl. Numer. Math. 44~(1) (2003) 139--181.

\bibitem{kennedy19:_higher_runge_kutta}
C.~A. Kennedy, M.~H. Carpenter, Higher-order additive {Runge-Kutta} schemes for ordinary differential equations, Appl. Numer. Math. 136 (2019) 183--205.

\bibitem{arndt21}
D.~Arndt, W.~Bangerth, D.~Davydov, T.~Heister, L.~Heltai, M.~Kronbichler, M.~Maier, J.-P. Pelteret, B.~Turcksin, D.~Wells, The \texttt{DEAL.II} finite element library: Design, features, and insights, Comput. Math. Appl. 81 (2021) 407--422.

\bibitem{burstedde11}
C.~Burstedde, L.~C. Wilcox, O.~Ghattas, \texttt{p4est}: Scalable algorithms for parallel adaptive mesh refinement on forests of octrees, SIAM J. Sci. Comput. 33~(3) (2011) 1103--1133.

\bibitem{kronbichler12}
M.~Kronbichler, K.~Kormann, A generic interface for parallel cell-based finite element operator application, Comput. \& Fluids 63 (2012) 135--147.

\bibitem{bangerth11:_algor}
W.~Bangerth, C.~Burstedde, T.~Heister, M.~Kronbichler, Algorithms and data structures for massively parallel generic adaptive finite element codes, ACM Trans. Math. Software 38~(2) (2011) 14:1--14:28.

\bibitem{clevenger21:_flexib_paral_adapt_geomet_multig_method_fem}
T.~C. Clevenger, T.~Heister, G.~Kanschat, M.~Kronbichler, A flexible, parallel, adaptive geometric multigrid method for {FEM}, ACM Trans. Math. Software 47~(1) (2020) 7:1--7:27.

\bibitem{taylor37:_mechan}
G.~I. Taylor, A.~E. Green, Mechanism of the production of small eddies from large ones, Proc. R. Soc. Lond. Ser. A, Math. Phys. Sci. 158~(895) (1937) 499--521.

\bibitem{nguyen11:_galer_navier_stokes}
N.~C. Nguyen, J.~Peraire, B.~Cockburn, An implicit high-order hybridizable discontinuous {Galerkin} method for the incompressible {Navier–Stokes} equations, J. Comput. Phys. 230~(4) (2011) 1147--1170.

\bibitem{ethier94:_exact_navier_stokes}
C.~R. Ethier, D.~A. Steinman, Exact fully {3D Navier-Stokes} solutions for benchmarking, Int. J. Numer. Methods Fluids 19~(5) (1994) 369--375.

\bibitem{piatkowski18:_galer_navier_stokes}
M.~Piatkowski, S.~Müthing, P.~Bastian, A stable and high-order accurate discontinuous {Galerkin} based splitting method for the incompressible {Navier–Stokes} equations, J. Comput. Phys. 356 (2018) 220--239.

\bibitem{doerfler96:_poiss}
W.~D\"orfler, A convergent adaptive algorithm for {P}oisson's equation, SIAM J. Numer. Anal. 33~(3) (1996) 1106--1124.

\bibitem{pfeiler20:_doerf}
C.-M. Pfeiler, D.~Praetorius, D\"orfler marking with minimal cardinality is a linear complexity problem, Math. Comp. 89~(326) (2020) 2735--2752.

\bibitem{bell91}
J.~B. Bell, P.~Colella, L.~H. Howell, An efficient second-order projection method for viscous incompressible flow, in: AIAA 10th Computational Fluid Dynamics Conference, Honolulu, Hawaii, June 24-26, 1991, pp. 360--367.

\bibitem{ghia82:_high_re_navier_stokes}
U.~Ghia, K.~N. Ghia, C.~T. Shin, {High-Re} solutions for incompressible flow using the {Navier-Stokes} equations and a multigrid method, J. Comput. Phys. 48~(3) (1982) 387--411.

\bibitem{erturk05:_numer_reynol}
E.~Erturk, T.~C. Corke, C.~G\"okc\"ol, Numerical solutions of {2-D} steady incompressible driven cavity flow at high {Reynolds} numbers, Int. J. Numer. Methods Fluids 48~(7) (2005) 747--774.

\bibitem{auteri02:_numer}
F.~Auteri, N.~Parolini, L.~Quartapelle, Numerical investigation on the stability of singular driven cavity flow, J. Comput. Phys. 183~(1) (2002) 1--25.

\bibitem{bruneau06}
C.-H. Bruneau, M.~Saad, The {2D} lid-driven cavity problem revisited, Comput. \& Fluids 35~(3) (2006) 326--348.

\bibitem{williamson96:_vortex_dynam_cylin_wake}
C.~H.~K. Williamson, Vortex dynamics in the cylinder wake, Annu. Rev. Fluid Mech. 28 (1996) 477--539.

\bibitem{braza86:_numer}
M.~Braza, P.~Chassaing, H.~H. Minh, Numerical study and physical analysis of the pressure and velocity fields in the near wake of a circular cylinder, J. Fluid Mech. 165 (1986) 79–130.

\bibitem{liu98:_precon_multig_method_unstead_incom_flows}
C.~Liu, X.~Zheng, C.~H. Sung, Preconditioned multigrid methods for unsteady incompressible flows, J. Comput. Phys. 139~(1) (1998) 35--57.

\bibitem{kolahdouz20}
E.~M. Kolahdouz, A.~P.~S. Bhalla, B.~A. Craven, B.~E. Griffith, An immersed interface method for discrete surfaces, J. Comput. Phys. 400 (2020) 108854.

\bibitem{calhoun02:_cartes_grid_method_solvin_two}
D.~Calhoun, A {Cartesian} grid method for solving the two-dimensional streamfunction-vorticity equations in irregular regions, J. Comput. Phys. 176~(2) (2002) 231--275.

\bibitem{mahir08:_numer}
N.~Mahír, Z.~Altaç, Numerical investigation of convective heat transfer in unsteady flow past two cylinders in tandem arrangements, Int. J. Heat Fluid Flow 29~(5) (2008) 1309--1318.

\bibitem{xu06}
S.~Xu, Z.~J. Wang, An immersed interface method for simulating the interaction of a fluid with moving boundaries, J. Comput. Phys. 216~(2) (2006) 454--493.

\bibitem{johnson99:_flow_reynol}
T.~A. Johnson, V.~C. Patel, Flow past a sphere up to a {Reynolds} number of 300, J. Fluid Mech. 378 (1999) 19--70.

\bibitem{ormieres99:_trans_turbul_wake_spher}
D.~Ormi\`eres, M.~Provansal, Transition to turbulence in the wake of a sphere, Phys. Rev. Lett. 83~(1) (1999) 80--83.

\bibitem{tomboulides00:_numer}
A.~G. Tomboulides, S.~A. Orszag, Numerical investigation of transitional and weak turbulent flow past a sphere, J. Fluid Mech. 416 (2000) 45--73.

\bibitem{hunt88:_eddies_stream_conver_zones_turbul_flows}
J.~C.~R. Hunt, A.~A. Wray, P.~Moin, Eddies, streams, and convergence zones in turbulent flows, {Proceedings} of the Summer Program, Center for Turbulence Research, Stanford University (1988).

\bibitem{bagchi00:_direc_numer_simul_flow_heat}
P.~Bagchi, M.~Y. Ha, S.~Balachandar, Direct numerical simulation of flow and heat transfer from a sphere in a uniform cross-flow, J. Fluids Eng. 123~(2) (2001) 347--358.

\bibitem{mittal08}
R.~Mittal, H.~Dong, M.~Bozkurttas, F.~M. Najjar, A.~Vargas, A.~{von Loebbecke}, A versatile sharp interface immersed boundary method for incompressible flows with complex boundaries, J. Comput. Phys. 227~(10) (2008) 4825--4852.

\bibitem{constantinescu03:_les_des_inves_turbul_flow}
G.~S. Constantinescu, K.~D. Squires, {LES} and {DES} investigations of turbulent flow over a sphere at {Re} = 10,000, Flow Turbul. Combust. 70~(1-4) (2003) 267--298.

\bibitem{kim01:_immer_bound_finit_volum_method}
J.~Kim, D.~Kim, H.~Choi, An immersed-boundary finite-volume method for simulations of flow in complex geometries, J. Comput. Phys. 171~(1) (2001) 132--150.

\bibitem{hartmann11:_cartes}
D.~Hartmann, M.~Meinke, W.~Schröder, A strictly conservative {Cartesian} cut-cell method for compressible viscous flows on adaptive grids, Comput. Methods Appl. Mech. Engrg. 200~(9) (2011) 1038--1052.

\bibitem{giacobello09:_wake_reynol}
M.~Giacobello, A.~Ooi, S.~Balachandar, Wake structure of a transversely rotating sphere at moderate {Reynolds} numbers, J. Fluid Mech. 621 (2009) 103--130.

\bibitem{burman23:_implic_oseen_reynol}
E.~Burman, D.~Garg, J.~Guzman, Implicit-explicit time discretization for {Oseen's} equation at high {Reynolds} number with application to fractional step methods, SIAM J. Numer. Anal. 61~(6) (2023) 2859--2886.

\bibitem{bansch16:_poster_error_estim_press_correc_schem}
E.~B\"{a}nsch, A.~Brenner, A posteriori error estimates for pressure-correction schemes, SIAM J. Numer. Anal. 54~(4) (2016) 2323--2358.

\bibitem{gross11:_numer_method_two_incom_flows}
S.~Gross, A.~Reusken, Numerical Methods for Two-phase Incompressible Flows, Springer-Verlag, Berlin Heidelberg, 2011.

\bibitem{richter17:_fluid_inter}
T.~Richter, Fluid-structure Interactions: Models, Analysis and Finite Elements, Springer, Cham, 2017.

\bibitem{zhang16:_mars}
Q.~Zhang, A.~Fogelson, {MARS}: An analytic framework of interface tracking via mapping and adjusting regular semi-algebraic sets, SIAM J. Numer. Anal. 54~(2) (2016) 530--560.

\bibitem{zhang18:_cubic_mars_method_fourt_sixth}
Q.~Zhang, Fourth- and higher-order interface tracking via mapping and adjusting regular semianalytic sets represented by cubic splines, SIAM J. Sci. Comput. 40~(6) (2018) A3755--A3788.

\bibitem{burman15:_cutfem}
E.~Burman, S.~Claus, P.~Hansbo, M.~G. Larson, A.~Massing, Cut{FEM}: Discretizing geometry and partial differential equations, Int. J. Numer. Meth. Engng. 104~(7) (2015) 472--501.

\bibitem{lehrenfeld19:_eulerian_fem}
C.~Lehrenfeld, M.~Olshanskii, An {E}ulerian finite element method for {PDE}s in time-dependent domains, ESAIM Math. Model. Numer. Anal. 53~(2) (2019) 585--614.

\end{thebibliography}

\end{document}